\documentclass[a4paper]{amsart}

 \usepackage{amsfonts}
  \usepackage{latexsym}
   \RequirePackage{amsbsy} 
   \RequirePackage{amsopn} 
   \RequirePackage{amsmath}
    \RequirePackage{amssymb}
    \usepackage{amsthm}

 \usepackage{color}%
 \DeclareMathOperator{\Spec}{Spec}
 \DeclareMathOperator{\Max}{Max}
 \DeclareMathOperator{\Inv}{Inv}
 \DeclareMathOperator{\Prin}{Prin}
 \DeclareMathOperator{\Cl}{Cl}
 \DeclareMathOperator{\Pic}{Pic}
  \DeclareMathOperator{\Ker}{Ker}
    \DeclareMathOperator{\Ima}{Im}
     \DeclareMathOperator{\G}{G}

\DeclareMathOperator{\calP} {\mathcal P}

\DeclareMathOperator{\calM} {\mathcal M}

   \newtheorem{thee}{Theorem}[section]
   \newtheorem{coor}[thee]{Corollary}
   \newtheorem{leem}[thee]{Lemma}

   \newtheorem{exxe}[thee]{Example}
   \newtheorem{reem}[thee]{Remark}
  \DeclareMathOperator{\oF}       {\boldsymbol{\overline{F}}}%
  \DeclareMathOperator{\F}        {\boldsymbol{F}}%

   \newcommand{\balf}
   {\renewcommand{\theenumi}{(\alph{enumi})}
   \renewcommand{\labelenumi}{\theenumi}
                        \begin{enumerate}}
  \newcommand{\ealf}   {\end{enumerate}
                        \renewcommand{\theenumi}{\arabic{enumi}}
                        \renewcommand{\labelenumi}{\theenumi.}}
  \newcommand{\bara}   {\renewcommand{\theenumi}{(\arabic{enumi})}
                        \renewcommand{\labelenumi}{\theenumi}
                        \begin{enumerate} }
  \newcommand{\eara}   {\end{enumerate}
                        \renewcommand{\theenumi}{\arabic{enumi}}
                        \renewcommand{\labelenumi}{\theenumi.}}

   \newcommand{\brom}   {\renewcommand{\theenumi}{(\roman{enumi})}
                        \renewcommand{\labelenumi}{\theenumi}
                        \begin{enumerate} }
  \newcommand{\erom}   {\end{enumerate}
                        \renewcommand{\theenumi}{\arabic{enumi}}
                        \renewcommand{\labelenumi}{\theenumi.}}

\begin{document}

\title{ON THE STAR CLASS GROUP OF A PULLBACK}

   \author{Marco Fontana  \;\;\;\;  \; Mi Hee Park}

   \dedicatory{Dipartimento di Matematica \\
   Universit\`a degli Studi Roma Tre \\
   {\rm \small Largo San Leonardo Murialdo, 1 }\\
   {\rm \small 00146 Roma, Italy} \\
   {\rm \small fontana@mat.uniroma3.it} \\
    \vspace{3pt}
    Department of Mathematics\\
    Chung-Ang University\\
    {\rm \small Seoul 156-756, Korea}\\
    {\rm \small mhpark@cau.ac.kr} }

 \date{\today}

 \thanks{\it Acknowledgment. \rm  During the preparation of this paper, the first named author
 was supported in part by a research
 grant MIUR 2003/2004.  }%

  \thanks{\it 2000 Mathematics Subject Classification. \rm   13C20, 13A15, 13G05.}%

\keywords{ Class group, Picard group, star operation, pullback, $t$--ideal, Pr\"ufer
multiplication domain.}%

\begin{abstract}
For the domain $R$ arising from the construction $T, M,D$,
we relate the star class groups of $R$ to those of $T$ and $D$.
More precisely, let $T$ be an integral domain, $M$ a nonzero maximal ideal of $T$,
$D$ a proper subring of $k:=T/M$, $\varphi: T\rightarrow k$ the natural projection,
and let $R={\varphi}^{-1}(D)$.
For each star operation $\ast$  on $R$, we define the star operation
$\ast_\varphi$ on $D$,  i.e.,  the ``projection''  of $\ast$ under $\varphi$, and
the star operation ${(\ast)}_{_{\!T}}$ on $T$,  i.e.,   the ``extension'' of $\ast$ to
$T$.
Then we show that, under a mild hypothesis on the group of units of
$T$, if $\ast$ is a star operation of finite type,
then the sequence of canonical homomorphisms
$0\rightarrow \Cl^{\ast_{\varphi}}(D) \rightarrow \Cl^\ast(R)
\rightarrow \Cl^{{(\ast)}_{_{\!T}}}(T)\rightarrow 0$
is split exact.  In particular, when $\ast =
t_{R}$, we deduce
that the sequence
$
0\rightarrow \Cl^{t_{D}}(D)
{\rightarrow} \Cl^{t_{R}}(R)
{\rightarrow}\Cl^{(t_{R})_{_{\!T}}}(T) \rightarrow 0
$
is split exact.  The relation between ${(t_{R})_{_{\!T}}}$ and
$t_{T}$ (and between $\Cl^{(t_{R})_{_{\!T}}}(T)$ and $\Cl^{t_{T}}(T)$)
is also investigated.
\end{abstract}

  \maketitle

\section{ Introduction and background results}

  The interest for constructing a general theory of the class group,
extending the theory of the  divisor  class group of a Krull domain,
was implicitly present already in the work by Claborn and Fossum (cf.
Fossum's book \cite{fo}).  One of the main objectives for this type of
extension was  to establish
a general functorial theory  by  exploiting class-group-type techniques in a more
general setting than that of Krull domains.
An approach to this problem, using star operations,
was initiated by D.F. Anderson in 1988 \cite{a}, where he studied in a
systematic way the star class group $\Cl^\star(R)$ of an integral
domain $R$, equipped with a star operation $\star$.  The key point of this
construction is that, when $\star$  is
the identity operation $d$,  $\Cl^d(R)$  coincides with the
Picard
group $\Pic(R)$ (which is, in fact, the ``classical'' class group of
the nonzero
fractional ideals when $R$ is a Dedekind domain);  when $\star$  is
the $v$--operation on a Krull domain,   $\Cl^v(R)$  coincides with
the ``usual''  divisor  class group of $R$;
when $\star$  is
the $t$--operation,   $\Cl^t(R)$, which is defined  on  arbitrary
domain $R$,  is commonly considered the best generalization of
the ``usual''  divisor  class group to the general setting (cf. the
pioneering work in this area by
Bouvier and Zafrullah \cite{bou82},   \cite{zafrullah85}, \cite{bz88}
and the recent excellent survey paper by D. F. Anderson \cite{a2}).

Since various divisibility properties are often reflected  in
group-theoretic properties of the class  groups,   a particular interest was
given in recent years to  the computation of the $t$--class
group  where  the functorial properties can be applied in
a very effective way (for instance, cf.  \cite{aaz}, \cite{gr}
and  \cite{nea}).

In  case  of  the   rings arising from  pullback
 construction  of various type (cf. \cite{f}, \cite{c}),  the $t$--class group  was
extensively studied by several authors
(cf.  for instance \cite{ar},
\cite{fg}, \cite{knea}, \cite {aebk}, \cite{eb},  and \cite {ac}).

It is  well known  that, even in the   case  of an embedding $A
\subset B$ of Krull domains, it is not possible in general to define a
canonical homomorphism
 between  the  divisor  class groups $\Cl(A) \rightarrow \Cl(B)$
(the condition (PDE), i.e.,  ``pas d'\'eclatement'', was introduced in
1964
by Samuel \cite{s} in order to characterize the existence of  this
canonical
homomorphism).
In case of star class groups,
 the technical difficulties for establishing functorial
properties were surmounted by D. F.
Anderson by introducing the notion of compatibility between star
operations. More precisely,  let $A$ be a subdomain of an integral
domain $B$ and let $\star_{_{\!A}} $ [respectively, $\star_{_{\!B}} $]
be a star operation on $A$ [respectively, on $B$], then  $\star_{_{\!A}} $
and $\star_{_{\!B}} $ are compatible  if
$(IB)^{\star_{_{\!B}}} = (I^{\star_{_{\!A}}}B)^{\star_{_{\!B}}}$ for each
nonzero fractional ideal $I$ of $A$.  In this situation,  the
extension map $I \mapsto IB$ induces a natural group homomorphism
$\Cl^{\star_{_{\!A}}}(A) \rightarrow \Cl^{\star_{_{\!B}}}(B)$.
Unfortunately, the compatibility condition is a sufficient but not
a necessary condition for the existence of the natural homomorphism
$\Cl^{\star_{_{\!A}}}(A) \rightarrow \Cl^{\star_{_{\!B}}}(B)$
\cite[page 823]{a}. Moreover,  the identity operation
$d_{_{\!A}}$ on $A$ is compatible with any star operation on $B$ while
it is very common that
the $t$--operation,  $t_{_{\!A}}$,  [respectively, the $v$--operation,
$v_{_{\!A}}$,] on $A$ is not compatible with
the $t$--operation, $t_{_{\!B}}$,  [respectively, the $v$--operation,
$v_{_{\!B}}$,] on $B$.
\smallskip

 In the present paper we mainly consider the following situation:

  \bf $(\square)$ \sl $\,T\,$ represents an integral domain,\ $M$\ a nonzero
  maximal ideal of \ $T$,\ $k$\ the residue field $T/M$,\ $D$\
      a  proper  subring  of\ $k$\ and $\varphi: T \rightarrow k$ the canonical
  projection.  Let $R: = \varphi^{-1}(D)=: T\times_{k }D$ be the integral
  domain arising from the following pullback of canonical homomophisms:
  $$
\begin{array}{ccc}
R & \longrightarrow & D \\
 \Big\downarrow & & \Big\downarrow
\\ T & \stackrel {\varphi}\longrightarrow & k=T/M.
\end{array}
$$ \rm
It is easy to see that $ M=(R:T)$ is the conductor of the
embedding $\iota: R \hookrightarrow T$.  In this situation, we will say
that we are dealing with \it a pullback of type \bf $(\square)$ \rm and
we will still denote by $\varphi$ the restriction
$\varphi\big\vert_{R}$, giving rise to a canonical surjective
homomorphism from $R=\varphi^{-1}(D)$ onto $D$.

Let $L$ denote the field of quotients of $D$ ( and hence,
  $L\subseteq k$). \sl  If we assume, moreover, that $L = k$,  \rm   then we will
  say that we are dealing with \it a pullback of type \bf
  $(\square^{+})$.\rm

The main goal of this work is to establish functorial relations
among the star class groups of $R$, $D$, and $T$, by using the
theory that we have recently  developed  in \cite{fp} concerning the
``lifting'' and  the  ``projection'' of a star operation under a
surjective homomorphism of integral domains, the ``extension'' of a
star operation to its overrings and the ``glueing'' of star
operations in pullback diagrams of a rather general type.   One of
the principal results proven  in this paper is that, given a
pullback
    diagram of type $(\square^{+})$ and a star operation $\ast$  of finite type  on
    $R$, if $\ast_{\varphi}$ denotes the ``projection''   of $\ast$ onto $D$
    [respectively, ${(\ast)}_{_{\!T}}$ denotes the ``extension'' of $\ast$  to
    $T$], under a mild hypothesis on the group of units of $T$,
the sequence of canonical homomorphisms
$$ 0\rightarrow \Cl^{\ast_{\varphi}}(D) \stackrel
{\overline{\boldsymbol{\alpha}}} {\longrightarrow} \Cl^\ast(R) \stackrel
{\overline{\boldsymbol{\beta}}}
{\longrightarrow}\Cl^{{(\ast)}_{_{\!T}}}(T)\rightarrow 0$$
is split exact (Theorem \ref{th:split}).  In particular, when $\ast =
t_{R}$, we deduce
that the sequence
$$
0\longrightarrow \Cl^{t_{D}}(D) \stackrel {\overline{\boldsymbol{\alpha}}}
{\longrightarrow} \Cl^{t_{R}}(R) \stackrel {\overline{\boldsymbol{\beta}}}
{\longrightarrow}\Cl^{(t_{R})_{_{\!T}}}(T) \longrightarrow 0
$$
is split exact.  The relation between ${(t_{R})_{_{\!T}}}$ and
$t_{T}$ (and between $\Cl^{(t_{R})_{_{\!T}}}(T)$ and $\Cl^{t_{T}}(T)$)
is also investigated.   Among the applications of the main results of
this paper, a characterization of when $R$ is a Pr\"ufer
$\ast$--multiplication domain is given.

\vskip.3in

  \centerline{ * * * * * }

  \vskip.3in \rm

    Let $D$ be an integral domain with quotient field $L$.  Let
   $\boldsymbol{\overline{F}}(D)$ denote the set of all nonzero
   $D$-submodules of $L$ and let $\boldsymbol{F}(D)$ be the set of all
   nonzero fractional ideals of $D$, i.e., all $E \in
   \boldsymbol{\overline{F}}(D) $ such that there exists a nonzero $d \in
   D$ with $dE \subseteq D$.  Let $\boldsymbol{f}(D)$ be the set of all
   nonzero finitely generated $D$-submodules of $L$.  Then, obviously $\,
   \boldsymbol{f}(D) \subseteq \boldsymbol{F}(D) \subseteq
   \boldsymbol{\overline{F}}(D) \, .$

   For each pair of fractional ideals $E, F$ of $D$, we denote as usual by
  $(E:_{L}F)$ the fractional ideal of $D$ given by $\{ y\in L\mid \, yF
  \subseteq E\}$; in particular, for each fractional ideal $I$
  of $D$, we set $I^{-1} := (D:_{L} I)$.
  \medskip

  We recall that a mapping
  $ \star : \boldsymbol{\overline{F}}(D) \rightarrow
  \boldsymbol{\overline{F}}(D) \,, \; E \mapsto E^{\star} $,
   is called a \it semistar operation on $D$ \rm if the following
   properties hold
   for all $0\neq x \in L$, and $E,F \in \boldsymbol{\overline{F}}(D)$:

  \hspace*{10pt}$(\star_1)$ \;\;  $(xE)^{\star} = xE^{\star}\,; $

  \hspace*{10pt}$(\star_2)$ \;\;  $E \subseteq F
  \Rightarrow  E^{\star} \subseteq F^{\star}\,; $

  \hspace*{10pt}$(\star_3)$\;\;  $E \subseteq E^{\star}$ and
  $E^{\star} =
  (E^{\star})^{\star} =: E^{\star \star} $

  \noindent (cf.  for instance \cite{OM1}, \cite{OM2},
  \cite{MSu}, \cite{MSa}, and  \cite{fh}).

  \begin{exxe} \label{ex:1.1}
  \bf (a) \rm If $\star$ is a semistar operation on $D$ such that $D^\star =
  D\,$, then the map (still denoted by) $\star : \boldsymbol{F}(D)
  \rightarrow \boldsymbol{F}(D)\,$, \ $E \mapsto E^\star ,$ is called a
  \it star operation on \rm $D\,.$ Recall \cite[(32.1)]{g} that a star
  operation $\star$ satisfies the properties $(\star_2)\,, (\star_3)\,$
  for all $E, F \in \boldsymbol{F}(D)\,;$ moreover, for each $0\neq x \in
  L$  and for each $E \in \boldsymbol{F}(D)\,,$ a star
  operation $\star$ satisfies the following ``stronger''  version of $(\star_1)$  (when restricted to $\boldsymbol{F}(D)$):

  \hspace*{10pt} $(\star\star_1)$ $(xD)^\star = xD\,; \; \; (xE)^{\star} = xE^{\star}\,.$

  Conversely, if $\star : \boldsymbol{F}(D)
  \rightarrow \boldsymbol{F}(D)\,$, \ $E \mapsto E^\star ,$ is a star
  operation on $D$ (i.e., if $\star$ satisfies the properties
  $(\star\star_1), \ (\star_2)\ $ and \  $(\star_3)$), then $\star$
  can be extended trivially to a semistar operation on $D$, denoted by $\star_{e}$ (or,
  sometimes, just by $\star$),  by
  setting $E^{\star_{e}} := L$, when $E \in
  \boldsymbol{\overline{F}}(D)\smallsetminus \boldsymbol{F}(D)$, and
  $E^{\star_{e}} := E^\star$, when $E \in
   \boldsymbol{F}(D)$.

 A semistar operation $\star$ on $D$ such that $D \subsetneq
  D^\star$ is called a \it proper semistar o\-pe\-ration on $D$.

  \bf (b) \rm The \it
  identity semistar operation $d_{D}$ on $\,D\,$ \rm (simply denoted by
  $d$) is a trivial semistar (in fact, star) operation on $D$ defined by
  $\, E^{d_{D}}: = E\,$ for each $\, E \in
  \boldsymbol{\overline{F}}(D)\,$ ($d_{D}$, when restricted to
  $\boldsymbol{{F}}(D)$, is a star operation on $D$).

   \bf (c) \rm For each $E \in \boldsymbol{\overline{F}}(D)$, set
  $E^{\star_{_{\!f}}} := \cup \{F^{\star} \;|\; \, F\subseteq E, \; F \in
  \boldsymbol{f}(D) \}\,.  $ \ Then $\star_{_{\!f}}$ is also a semistar operation
  on $D$, which is called \it the semistar operation of finite type
  associated to $\star\,$.  \rm Obviously, $F^{\star} = F^{\star_{_{\!f}}}$ for
  each $F \in \boldsymbol{f}(D)\, ;$ \ moreover, if $\star$ is a star
  operation, then ${\star_{_{\!f}}}$ is also a star operation.  If $\star =
  \star_{_{\!f}}$, then the semistar [respectively, the star] operation $\star$
  is called a \it semistar \rm [respectively, \it star\rm ] \it operation
  of finite type.  \rm

   Note that $\star_{_{\!f}}\leq \star\,,$  i.e., $\, E^{\star_{_{\!f}}}
  \subseteq E^{\star}\, $ for each $E \in \boldsymbol{\overline{F}}(D)$.
  Thus, in particular, if $E = E^{\star}$, then $E = E^{\star_{_{\!f}}}$.  Note
  also that $\star_{_{\!f}} = (\star_{_{\!f}})_{_{\!f}}$.

  More generally, if ${\star_1}$ and ${\star_2}$ are two semistar operations on $D$,
  we say that \ ${\star_1} \leq {\star_2}$ \ if $E^{\star_1} \subseteq
  E^{\star_2}$ for each $E \in \boldsymbol{\overline{F}}(D)$.  In this situation,
  it is easy to see that $\left(E^{\star_{1}}\right)^{\star_{2}} =
  E^{\star_2}= \left(E^{\star_{2}}\right)^{\star_{1}}$.

  There are several examples of nontrivial semistar or star operations of
  finite type; the best known is probably the $t$--operation.  Indeed, we
  start from the \it $v_{D}$ star operation \rm on an integral domain $D$
  (simply denoted by $v$), which is defined by

  \centerline{ $E^{v_{D}} := (E^{-1})^{-1}=(D:_{L}(D:_{L}E))\, $}

  \noindent for any $E\in\F(D)$, and we set  $t_{D}:=(v_{D})_{_{\!f}}$ (or, simply, $t
  =v_{_{\!f}}$).

  \bf (d) \rm Let $\iota:R \hookrightarrow T$ be an embedding of integral
  domains with the same field of quotients $K$ and let $\ast$ be a
  semistar operation on $R$.  Define $\ast_{\iota} :
  \overline{\boldsymbol{F}}(T) \rightarrow \overline{\boldsymbol{F}}(T)$
  by setting

  \centerline{$ E^{\ast_{\iota}} := E^\ast\, \; \mbox{ for each } \;
   E \in \overline{\boldsymbol{F}}(T) \
  (\subseteq \overline{\boldsymbol{F}}(R) )\,.  $}

  \noindent Then,  it is easy to verify (cf. also \cite[Proposition
  2.8]{FL1}) that:

   \hskip 0.5cm  \bf (d1) \sl If $\iota$ is not the identity map, then
   $\ast_{\iota}$ is a semistar, possibly non--star, operation on $T$, even
   if $\ast$ is a star operation on $R$ (obviously, if $\iota$ is the identity map, then $\ast_{\iota} = \ast$ and thus this phenomenon does not occur) .  \rm

     \rm Note that when $\ast$ is a star operation on $R$ and
     $(R:_{K}T) =(0)$, a fractional ideal $E$
     of $T$ is  not  a fractional ideal of $R$, hence $\ast_{\iota}$
     is not  necessarily  defined as a star operation on $T$.

  \hskip 0.5cm  \bf (d2) \sl If $\ast$ is of finite type
  on $R$, then
  ${\ast_{\iota}} $ is also of finite type    on $T\,.$

  \hskip 0.5cm  \bf (d3) \sl If $T := R^\ast$, then ${\ast_{\iota}}$ defines
  a star operation on      $T$.

    \bf (e) \rm Let $\star$ be a semistar operation on the overring $T$ of
    $R$.  Define $\star^{\iota}: \overline{\boldsymbol{F}}(R) \rightarrow
    \overline{\boldsymbol{F}}(R)$ by setting

  \centerline{$ E^{\star^{\iota}} := (ET)^\star\, \; \;\; \; \mbox{ for each } \; E \in
  \overline{\boldsymbol{F}}(R)\,.
  $}

  \noindent Then, we know \cite[Proposition 2.9, Corollary 2.10]{FL1}:

  \hskip 0.5cm \bf (e1) \sl ${\star^{\iota}} $ is a semistar operation on $R$.

   \hskip 0.5cm \bf (e2) \sl If $\star:= d_{T}$, then ${(d_{T})}^{\iota}$ is a
   semistar operation of finite type on $R$, \rm     which is also denoted by
   $\star_{\{T\}}$ (i.e., it is the semistar operation on $R$ defined by
   $E^{\star_{\{T\}}} := ET$ for each $E \in \oF(R)$).

  \hskip 0.5cm \bf (e3) \sl For each semistar operation $\star$ on $T$,
$({\star^{\iota}})_{\iota}= \star$.

  \bf (f) \rm Let $\Delta $ be a set of prime ideals of an
  integral domain $D$ with quotient field $L$.  The mapping $E \mapsto E^{\star_\Delta}$,
  where $E^{\star_{\Delta}}:=\cap \{ ED_P\; | \;\,
  P\in \Delta \}$\, for each $E\in \oF(D)$,
  defines a semistar operation on $D\,.$ Note that ${\star_{\Delta}}$
  (restricted to the nonzero fractional ideals of $D$) is a star
  operation on $D$ if and only if
  $D=\cap \{ D_P\; | \;\, P\in \Delta \}$.  Moreover (\cite [Lemma
  4.1]{fh} or \cite [Theorem 1]{dda}):

  \hskip 0.5cm \bf (f1) \sl For each $E\in \oF(D)$ and for each $P\in
  \Delta\,$, $ED_P=E^{\star_\Delta}D_P$.

  \hskip 0.5cm \bf (f2) \sl The semistar operation $\star_\Delta$ is \rm
  stable (with respect to the finite intersections), \it i.e., for all $E,
  F \in \oF(D)\,,$ we have $\, (E \cap F)^{\star_{\Delta}} =
  E^{\star_{\Delta}} \cap F^{\star_{\Delta}}\,.  $

  A \it semistar operation $\star$ on $D$ \rm is called \it spectral \rm if
  there exists a subset $\Delta $ of Spec($D$) such that $\star =
  \star_{\Delta}\,;$\ in this case we say that $\star$ is \it the spectral
semistar operation associated with \rm $\Delta \,.$

\bf (g) \rm Let $\star$ be a star operation on $D$.  If $E \in
\boldsymbol{{F}}(D)$, we say that $E$ is \it a $\star$--ideal \rm if
$E=E^\star$.  We denote by $\boldsymbol{{F}}^{\star}(D)$ (respectively,
$\boldsymbol{{f}}^{\star}(D)$) the set $\{E\in \boldsymbol{{F}}(D) \mid
E=E^\star\}$ (respectively, $\{E\in \boldsymbol{{F}}(D) \mid
E=F^\star\,\ \mbox {where } F \in \boldsymbol{{f}}(D)\}$.  Obviously,
$\boldsymbol{{F}}^d(D)=\boldsymbol{{F}}(D)$ (respectively,
$\boldsymbol{{f}}^d(D)=\boldsymbol{{f}}(D)$) and the set
$\boldsymbol{{F}}^v(D)$ is called \it the set of divisorial ideals of
$D$.\rm

  Set $\calP(\star) := \Spec^\star(D):= \{P\in \Spec(D) \mid P =
  P^\star\}$  and $\calM(\star) := \Max^\star(D)$ which is
  the     (possibly empty) set      of all the maximal
  elements of the set $\{ I $ proper ideal of $ D \mid I =
  I^{\star} \}$.
  Assume that each proper $\star$--ideal of $D$ is contained in some prime
  ideal of $\Spec^\star(D)$, then it is known that $\star_{\calP(\star)}$ is a
     star operation      on $D$ \cite[Theorem 3]{dda}.  In particular, for each star
  operation $\star$  on $D$ which is not a field,
   $\calM(\star_{_{\!f}})$  is a nonempty subset of $\calP(\star_{_{\!f}})$
  and it satisfies the property that each proper $\star_{_{\!f}}$--ideal of
  $D$ is contained in some prime ideal of $\calM(\star_{_{\!f}})$.  Then
  $\widetilde{\star} := \star _{\calM(\star_{_{\!f}})}$ is a star
  operation of finite type and stable on $D$, which is called \it the stable
  operation of finite type associated to $\star$.  \rm It is easy to see that
  $\widetilde{\ \star_{_{\!f}}} =\widetilde{\star} =
  (\widetilde{\star})_{_{\!f}}$ and $\widetilde{\star} =
  \star_{\calP(\star_{_{\!f}})}$.  Note that \cite[Corollary 3.9]{fh}

  \centerline{ $\star = \widetilde{\star} \; \Leftrightarrow \; \star
  \mbox{ \sl     is a stable star operation of finite type.     \rm}$}

Particularly interesting is the case in which $\star =v$.
Using the notation introduced by Wang Fanggui and
  R.L. McCasland \cite{Fanggui/McCasland:1997}, we will denote by $w_{D}$
  (or, simply, $w$) the star operation $\widetilde{v_{D}} =
  \widetilde{t_{D}}$ (simply, $w := \widetilde{v}=
  \widetilde{t}$;  cf.  also \cite{Hedstrom/Houston: 1980} and
  \cite{ack}).

  Note that if ${\star_1}$ and ${\star_2}$ are two star operations on $D$,
  then
$$
{\star_1} \leq
{\star_2} \;\; \Leftrightarrow \;\; \boldsymbol{{F}}^{\star_{2}}(D)
\subseteq \boldsymbol{{F}}^{\star_{1}}(D)\,.
$$
It is well known that for each star operation $\star$, we have
$\widetilde{\star} \leq \star_{\!_f} \leq \star$ \cite[Theorem
2.3]{ack}.  Thus, in particular, if $E = E^{\star}$, then $E =
E^{\widetilde{\star}} =E^{\star_{_{\!f}}}$.  Moreover, note that
$$
\boldsymbol{{f}}^{\star}(D) = \boldsymbol{{f}}^{\star_{_{\!f}}}(D) \subseteq
\boldsymbol{{F}}^{\star}(D) \subseteq \boldsymbol{{F}}^{\star_{_{\!f}}}(D)\,.
$$
It is also known that if ${\star_1}$ and ${\star_2}$ are two star
operations on $D$ and ${\star_1} \leq {\star_2}$, then
$({\star_1})_{_{\!f}}\leq ({\star_2})_{_{\!f}}$ and $\widetilde{{\
\star_1}}\leq \widetilde{{\ \star_2}}$.  In particular, for each
star operation $\star$, we have $\star \leq v$ \cite[Theorem 34.1
(4)]{g} and so $\star_{\!_f} \leq t$ and $\widetilde{\star} \leq w$.
Thus we get
$$\begin{array}{rl} \boldsymbol{{F}}^{v}(D) \subseteq&\hskip -6pt
\boldsymbol{{F}}^{t}(D) \subseteq\boldsymbol{{F}}^{w}(D) \subseteq
\boldsymbol{{F}}(D)\,,\\
\boldsymbol{{F}}^{v}(D) \subseteq&\hskip -6pt \boldsymbol{{F}}^{\star}(D)\,,\;\;\;
 \boldsymbol{{F}}^{t}(D) \subseteq \boldsymbol{{F}}^{\star_{_{\!f}}}(D)\,,\;\;\;
 \boldsymbol{{F}}^{w}(D)
\subseteq  \boldsymbol{{F}}^{\widetilde{\star}}(D)\,.
\end{array}
$$

  \bf (h) \rm Let $\iota:R \hookrightarrow T$ be an embedding of integral
  domains with the same field of quotients $K$ and let $\ast$ be a
  semistar operation on $R$.    It is not difficult to prove:

   \centerline{ \sl  $\ast$ is
  stable on $R$  $\; \Rightarrow \; $
  ${\ast_{\iota}} $ is stable on $T\,.$} \rm

\bf (k) \rm \rm If $\{\star_{\lambda}\mid \lambda \in \Lambda\}$ is a
family of semistar [respectively, star] o\-pe\-rations on $D$, then \
$\wedge_{\lambda} \{\star_{\lambda}\mid \lambda \in \Lambda\}$ \ (    simply denoted
by      \ $\wedge \star_{\lambda}$\ ), defined by

       \centerline{
       $ E^{\wedge \star_{\lambda}}:= \cap \{ E^{\star_{\lambda}}\mid \lambda
       \in \Lambda\}\,,$ \;\;\;\ {\mbox{for each $E \in
       \overline{\boldsymbol{F}}(D)$ \ [respectively, $E \in
       \boldsymbol{F}(D)$]\,,} } }

    \noindent is a semistar [respectively, star] operation on $D$.  Note
    that if at least one of the semistar operations in the family
    $\{\star_{\lambda}\mid \lambda \in \Lambda\}$ is a star operation on $D$,
    then\ $\wedge \star_{\lambda}$\ is still a star operation on $D$.\rm
\end{exxe}

  \vskip.1in  \rm

  Let $\star$ be a star operation on an integral domain $D$ and let $F \in
  \boldsymbol{F}(D)$. We say that $F$ is \emph{$\star$--invertible} if
  $\left(FF^{-1}\right)^\star=D$.  In particular, when $\star = d$
  [respectively,\ $v$\,,\; $t$\,,\; $w$\ ] is the identity star operation
  [respectively,\ the $v$--operation,\; the $t$--operation,\; the
  $w$--operation\ ], we reobtain the classical notion of \it invertibility
  \rm [respectively,\ \it $v$--invertibility,\; $t$--invertibility,\;
  $w$--invertibility\ \rm] of a fractional ideal.  Recall that:

\begin{leem} \label{lemma:inv}
Let $\star, {\star_1}, {\star_2}$ be     star      operations on an
integral domain $D$. Let  $\Inv(D, \star)$ be the set of all
$\star$--invertible fractional ideals of $D$ and $\Inv(D)$ (instead of
$\Inv(D, d)$) the set of all invertible fractional ideals of $D$.  Then
\begin{enumerate}
\item[(1)] $D \in  \Inv(D, \star)$.

\item[(2)] If ${\star_1} \leq {\star_2}$, then $ \Inv(D, \star_1) \subseteq
\Inv(D, \star_2)$. In particular, $\Inv(D) \subseteq \Inv(D,
\tilde{\star}) \subseteq \Inv(D,\star_{\!_f}) \subseteq \Inv(D,
\star)$ and so $\Inv(D) \subseteq \Inv(D,
w) \subseteq \Inv(D, t) \subseteq \Inv(D, v)$.

\item[(3)] $I,J \in  \Inv(D, \star)$ if and only if $IJ \in
 \Inv(D, \star)$.

\item[(4)] If $I \in  \Inv(D, \star)$, then $I^{-1} \in
 \Inv(D, \star)$.

\item[(5)] If $I \in  \Inv(D, \star)$, then $I^v \in
 \Inv(D, \star)$.
\end{enumerate}
\end{leem}

Let $\star$ be a star operation on $D$. Then $\boldsymbol{F}^{\star}(D)$ is a commutative monoid
under the $\star$--multiplication defined by $(I,J)\mapsto (IJ)^\star\,$
for each $I,\,J\in \boldsymbol{F}^{\star}(D)$.
If $\star_1$ and $\star_2$ are two star operations on $D$ with $\star_1\leq \star_2$, then
while $\boldsymbol{F}^{\star_2}(D)\subseteq \boldsymbol{F}^{\star_1}(D)$,
$\boldsymbol{F}^{\star_2}(D)$ is not a submonoid of $\boldsymbol{F}^{\star_1}(D)$ in general
(see \cite[page 811]{a}).
However, there is a special submonoid of $\boldsymbol{F}^{\star}(D)$ which reverses the inclusion:

\begin{leem} \label{lemma:1.2} \rm  (D.F. Anderson \cite[Proposition
3.3]{a}). \it Let $\star, {\star_1}, {\star_2}$ be star operations on an
integral domain $D$ and suppose that ${\star_1}\leq {\star_2}$.  Let
$\Inv^\star(D) := \{I \in \Inv(D, \star) \mid I = I^\star\}$ be the
set of all $\star$--invertible $\star$--ideals of $D$ and let $\Inv(D)$
(instead of $\Inv^{d}(D)$) be the set of all invertible fractional ideals
of $D$.  Then
\begin{enumerate}
\item[(1)] $\Inv^\star(D)$ is a submonoid of $\boldsymbol{F}^{\star}(D)$; moreover, it is an abelian group.
\item[(2)]   $\Inv^{\star_1}(D)$ is a subgroup of $\Inv^{\star_2}(D)$ (in symbol, $\Inv^{\star_1}(D) \leq
\Inv^{\star_2}(D)$).     In particular, for each star operation $\star$
on $D$,\ $\Inv(D) \leq \Inv^\star(D) \leq \Inv^v(D)$,\ $\Inv(D) \leq
\Inv^{\star_{_{\!f}}}(D) \leq \Inv^t(D)$\ and\ $\Inv(D) \leq
\Inv^{\widetilde{\star}}(D) \leq \Inv^{\star_{_{\!f}}}(D) \leq
\Inv^{\star}(D)\,.$  
\end{enumerate}
\end{leem}

 \rm In \cite{fp} we considered the problem of ``lifting a star
 operation'' with respect to a surjective ring homomorphim between two
 integral domains.  More precisely:

  \begin{leem} \label{le:1.3} \rm \cite[Corollary 2.4]{fp}.\ \it Let $R$ be an
  integral domain with field of quotients $K$, $M$ a prime ideal of $R$.
  Let $D$ be the quotient-domain $R/M$ and let $\varphi: R \rightarrow D$
  be the canonical projection.  Assume that $\star$ is a star operation on
  $D$.  For each nonzero fractional ideal $E$ of $R$, we set $$\begin{array}
  {rl} E^{\star^{\varphi}} &:= \cap \left\{
  x^{-1}\varphi^{-1}\left(\left(\frac{xE+M}{M}\right)^{\star}\right) \mid
  \, x \in E^{-1}\,, \; x \neq 0 \right\}\\ &= \cap \left\{
  x\varphi^{-1}\left(\left(\frac{x^{-1}E+M}{M}\right)^{\star}\right) \mid
  \, x \in K\,, \; E \subseteq xR \right\}\,,
  \end{array}
  $$
  where, if\,  $\frac{zE+M}{M}$\,  is the zero ideal of $D$, then we set
  $\varphi^{-1}\left(\left(\frac{zE+M}{M}\right)^{\star}\right)$ = $M$.
  Then ${\star^{\varphi}}$ is a star operation on $R$.  
      \end{leem}

       \rm In \cite{fp} we also considered the problem of  ``projecting a star
       operation'' with respect to  a  surjective homomorphism of
       integral domains, with particular emphasis on pullback constructions of
       a ``special'' kind.  More precisely:

  \begin{leem} \label{le:1.4} \rm \cite[Propositions 2.6, 2.7, 2.9 and
  Theorem 2.12]{fp}.\  \it Let $\varphi: R \rightarrow D$ be a surjective
  homomorphism of integral domains, let $\ast$ be a star operation on
  $R$ and let $L$ be the quotient field of $D$.
  For each nonzero fractional ideal $F$ of $D$, we set

\centerline{$ F^{\ast_{\varphi}}:= \cap
  \left\{y\varphi\left(\left(\varphi^{-1}\left(y^{-1}F\right)\right)^{\ast}\right)
  \mid \, y \in     L      \,,\;  F \subseteq yD \right\} \,.$}
  \begin{enumerate}
\item[(1)] \ $\ast_{\varphi}$ is a star operation on $D$.

\hskip -1.25cm   Assume, now, that we are dealing with a pullback diagram of
type  $(\square)$.  Then

\item[(2)] \  $ F^{\ast_{\varphi}} =
\varphi\left(\left(\varphi^{-1}(F)\right)^\ast\right) =
\left(\varphi^{-1}(F)\right)^\ast/M\,$ for each $F \in \boldsymbol{F}(D)$.

\item[(3)] \
$(\star^\varphi)_{\varphi} = \star\,$ for each star operation $\star$ on $D$.

\item[(4)] \
$\ast \leq (\ast_{\varphi})^{\varphi}\,$ for each star operation $\ast$
on $R$.  
\end{enumerate}
  \end{leem}

\section{Main results}

\begin{leem} \label{le:2.1}
    Assume that we are dealing with
a pullback diagram of type $(\square^{+})$.  Let $\ast$ be a star
operation on $R$ and let $\ast_{\varphi}$ be the star operation on $D$
defined in\ Lemma \ref{le:1.4}.  Then the map $\boldsymbol{\alpha}(\varphi,
\ast) $ (or, simply, $\boldsymbol{\alpha}$)$: \Inv(D, \ast_{\varphi})
\rightarrow \Inv(R, \ast)$, defined by\ $J \mapsto \varphi^{-1}(J)$, is
injective with $\Ima(\boldsymbol{\alpha}) = \{I\in \Inv(R, \ast)\mid
M\subsetneq I \subseteq I^{v_R}\subsetneq T\}$.  Moreover, if we use the
same notation $\boldsymbol{\alpha} =\boldsymbol{\alpha}(\varphi, \ast)$
for the restriction of the map $\boldsymbol{\alpha}$ to the subset
$\Inv^{\ast_{\varphi}}(D)$, then $\boldsymbol{\alpha}:
\Inv^{\ast_{\varphi}}(D) \rightarrow \Inv^{\ast}(R) $ is still
injective with $\Ima(\boldsymbol{\alpha}) = \{I\in \Inv^{\ast}(R) \mid
M\subsetneq I \subseteq I^{v_R}\subsetneq T\}$.
\end{leem}

\begin{proof} Recall first that the map $J\mapsto \varphi^{-1}(J)$ establishes a 1-1 correspondence
between $\boldsymbol{F}(D)$ and the set $\{H\in \boldsymbol{F}(R)\mid
M\subsetneq H \subseteq H^{v_R}\subsetneq T\}$ \cite[Corollary 1.9]{fg}.
Let $J\in \boldsymbol{F}(D)$. Then by applying Lemma \ref{le:1.4} (2),
we have $J^{\ast_{\varphi}} = {\left(\varphi^{-1}(J)\right)^\ast}/{M}$.
Therefore, $$ J =
J^{\ast_{\varphi}} \;\; \Leftrightarrow\;\; \varphi^{-1}(J)=
\left(\varphi^{-1}(J)\right)^\ast\,,$$
$$
(JJ^{-1})^{\ast_{\varphi}} = D\; \Leftrightarrow\;
\left(\varphi^{-1}\left(JJ^{-1}\right)\right)^\ast =R \,.
$$
By \cite[Proposition 1.6 and Proposition 1.8 (a)]{fg},
$\varphi^{-1}(JJ^{-1})=\varphi^{-1}(J)\varphi^{-1}(J^{-1})
=\varphi^{-1}(J)\left(\varphi^{-1}(J)\right)^{-1}$.
Therefore, $$
(JJ^{-1})^{\ast_{\varphi}} = D\; \Leftrightarrow\;
\left(\varphi^{-1}\left(J\right)\left(\varphi^{-1}(J)\right)^{-1}\right)^\ast
= R\,.$$

\vskip -18pt\end{proof}

Let $\Prin(D)$ be the subgroup of $\Inv^\star(D)$
of all \it the nonzero fractional principal ideals of $D$. \rm We recall that the
quotient group
$$
\Cl^\star(D): = \frac{\Inv^\star(D)}{\Prin(D)}
$$
is called \it the class group of an integral domain $D$ with
respect to a star operation $\star$ on $D$.  \rm

If $\star =d$ is the
identity star operation on $D$, then $\Cl^d(D)$ is denoted by $\Pic(D)$
and it is called \it the Picard group of an integral domain $D$.  \rm

\begin{leem} \label{le:2.2}
    Let $\star, {\star_1}, {\star_2}$ be star operations on
an integral domain $D$ and suppose that ${\star_1}\leq {\star_2}$.
Then $\Cl^{\star_1}(D)$ is a subgroup of $\Cl^{\star_2}(D)$.\ In
particular, for each star operation $\star$ on $D$,\ $\Pic(D) \leq
\Cl^\star(D) \leq \Cl^v(D)$, $\Pic(D) \leq \Cl^{\star_{_{\!f}}}(D) \leq
\Cl^t(D)$ and $\Pic(D) \leq \Cl^{\widetilde{\star}}(D)
\leq\Cl^{\star_{_{\!f}}}(D) \leq \Cl^{\star}(D)$.
\end{leem}

\begin{proof} Easy consequence of Lemma \ref{lemma:1.2}.  \end{proof}

    \begin{reem} \label{rk:2.3} \rm Note that the previous statement can be
   strengthen, since Anderson-Cook (in \cite[Theorem 2.18]{ack}) proved that
   \sl for any star operation $\star$ on an integral domain $D$,
   $\Inv^{\widetilde{\star}}(D) =\Inv^{\star_{_{\!f}}}(D)$, and thus
   $\Cl^{\widetilde{\star}}(D)=\Cl^{\star_{_{\!f}}}(D)$.
    \end{reem} \rm

\begin{leem}\label{lemma2}
 Assume that we are dealing with
a pullback diagram of type $(\square^{+})$.  Then the following
statements are equivalent:
\begin{enumerate}
\item[(1)] the canonical map $\tilde{\varphi}: \mathcal U(T)\rightarrow
k^\bullet/\mathcal U(D)$, $u\mapsto \varphi(u)\mathcal U(D)$, is a surjective
group homomorphism, where $k^\bullet$ is the multiplicative group of the
nonzero elements of the field $k$ and  $\mathcal U(T)$ (respectively, $\mathcal
U(D)$) is the group of units of $T$ (respectively, $D$);

\item[(2)] for
each nonzero element $x\in k$, ${\varphi}^{-1}(xD)$ is a fractional
principal ideal of $R$;

\item[(3)] the map $\overline{\boldsymbol{\alpha}}(\varphi, \ast)$ (or, simply,
$\overline{\boldsymbol{\alpha}}$) : $\Cl^{\ast_{\varphi}}(D) \rightarrow \Cl^\ast(R)$,
$[J]\mapsto [{\varphi}^{-1}(J)]\ (= [\boldsymbol{\alpha}(J)]$, where
$\boldsymbol{\alpha}$ is defined in Lemma \ref{le:2.1}), is a
well-defined group homomorphism for any star operation $\ast$ on $R$.
\end{enumerate}
\end{leem}

\begin{proof} (1) $\Leftrightarrow$ (2) $\Leftarrow$ (3)
see \cite[Theorem 2.3 (i)$\Leftrightarrow$(ii)$\Leftarrow$(iv)]{fg}.
    The direction (2) $\Rightarrow$ (3) is a consequence of Lemma
    \ref{le:2.1}.
\end{proof}

\begin{reem} \rm General examples for which  the map $\tilde{\varphi}:
    \mathcal U(T)\rightarrow k^\bullet/\mathcal U(D)$ is surjective are provided in \cite[Proposition 2.9]{fg}.
    \end{reem} \rm

The next theorem provides a generalization of the result by D. F. Anderson
\cite[Proposition 5.5]{a}:

\begin{thee} \label{th:2.4} Assume that we are dealing with a pullback diagram of
type $(\square)$.  If, moreover, $T$ is quasilocal, then the canonical
map $\overline{\boldsymbol{\alpha}}\ ( =
\overline{\boldsymbol{\alpha}}(\varphi, \ast)) : \Cl^{\ast_{\varphi}}(D)
\rightarrow \Cl^\ast(R)$ is an isomorphism  for any star operation $\ast$
on $R$.
\end{thee}

\begin{proof} We adapt the argument used in the proof of \cite[Proposition
5.5]{a}.  We first show that  $\Cl^\ast(R)=0$ when $D$ is a proper subfield
of $k$.  In this case, $R$ is quasilocal, since
$R$ and $T$ have the same prime spectrum \cite{ad}.  Let $I\in
\Inv^\ast(R)$.  As $M=(R:T)$ is a divisorial ideal of $R$, if
$II^{-1}\subseteq M$, then $(II^{-1})^\ast\subseteq M^\ast=M$, a
contradiction.  Then, necessarily, $II^{-1}=R$; thus $I$ is invertible
in the quasilocal domain $R$, and hence $I$ is principal.  Thus
$\Cl^\ast(R)=0$.

Without loss of generality, we may assume that $D$ is a proper subring of
$k$ with quotient field $k$, i.e., that we are dealing with
a pullback diagram of type $(\square^{+})$.  In this situation,  the map
$\overline{\boldsymbol{\alpha}}$: $\Cl^{\ast_{\varphi}}(D) \rightarrow
\Cl^\ast(R)$ is a homomorphism, because when $T$ is quasilocal, the condition
(1) of Lemma ~\ref{lemma2} holds \cite[Proposition 2.9]{fg}.

Let $J\in
\Inv^{\ast_{\varphi}}(D)$ such that ${\varphi}^{-1}(J)$ is principal in
$R$, say ${\varphi}^{-1}(J)= xR\,$ for some nonzero $x\in T$.  Then
$J=xR/M=\varphi(x)D$ is principal in $D$.  Therefore
$\overline{\boldsymbol{\alpha}}$ is injective.

Conversely, let $I\in \Inv^\ast(R)$.  Then, necessarily,
$II^{-1}\not\subseteq M$, and hence $II^{-1}T=T$, i.e., $IT$ is
invertible in $T$.  Since $T=R_M$ is quasilocal \cite[Corollary 0.5]{fg},
$IT=IR_M$ is principal, say $IT=iR_M$ for some $i\in I$.  Set
$I_1:=i^{-1}I$.  Then, obviously, $I_1\in \Inv^\ast(R)$ and $R\subseteq
I_1\subseteq T= I_1T$.  To prove that $\varphi(I_{1})=I_1/M$ belongs to
$\Inv^{\ast_{\varphi}}(D)$, it suffices to show that $(I_1)^v\subsetneq
T$ by Lemma ~\ref{le:2.1},       because
${\varphi}^{-1}\left({\varphi}(I_{1})\right) =I_{1}$.
 Suppose  that
$(I_1)^v=T$, then $I_1^{-1}=(R:T)=M$.  So
$R=\left(I_1I_1^{-1}\right)^{\ast} = (I_1M)^\ast\subseteq
(TM)^\ast=M^\ast=M$, a contradiction.  Thus, necessarily, we have
$(I_1)^v\subsetneq T$.  Therefore
$[I]=[i^{-1}I]=[I_1]=[{\varphi}^{-1}(I_1/M)]=\overline{\boldsymbol{\alpha}}([I_1/M])$.
 Hence $\overline{\boldsymbol{\alpha}}$ is also surjective and thus we
conclude that $\overline{\boldsymbol{\alpha}}$ is an isomorphism.
\end{proof}

\begin{coor} \label{cor:2.5} Assume that we are dealing with a pullback diagram of
type $(\square)$.  If, moreover, $T$ is quasilocal, then we have the
following canonical isomorphisms: $$\Pic(D)\cong \Pic(R)\,, \;\;\;
\Cl^t(D)\cong \Cl^t(R)\,, \;\;\; \Cl^w(D)\cong \Cl^w(R)\,, \;\;\;
\Cl^v(D)\cong \Cl^v(R)\,.$$
\end{coor}
\begin{proof}  Since $(d_R)_{\varphi}=d_D$, $(t_R)_{\varphi}=t_D$,
$(w_R)_{\varphi}=w_D$ and $(v_R)_{\varphi}=v_D$ (\cite[Proposition
3.3, Proposition 3.7, Corollary 3.10, and Corollary 2.13]{fp}), the
conclusion follows from the above theorem. The third isomorphism
also follows from the second one  by  Remark \ref{rk:2.3}.
\end{proof}

\begin{coor} \label{cor:2.7} Assume that we are dealing with
a pullback diagram of type $(\square)$. Let $T$ be quasilocal. Then
\begin{enumerate}
\item[(1)] The canonical homomorphism
$\overline{\boldsymbol{\alpha}}(\varphi, {\star}^{\varphi}):
\Cl^{\star}(D) \rightarrow \Cl^{{\star}^{\varphi}}(R)$ is an
isomorphism for any star operation $\star$ on $D$.
\item[(2)] $\Cl^\ast(R) =\Cl^{({\ast}_{\varphi})^{\varphi}}(R)$
for any star operation $\ast$ on $R$.
\end{enumerate}
\end{coor}
\begin{proof} (1) Set $\ast:= {\star}^{\varphi}$.  Then
$\ast_{\varphi}=({\star}^{\varphi})_{\varphi}=\star$ by Lemma \ref{le:1.4} (3).  The
conclusion follows immediately from Theorem \ref{th:2.4}.

(2) Recall that $\ast \leq  (\ast_{\varphi})^{\varphi}$ and
$((\ast_{\varphi})^{\varphi})_{\varphi}=\ast_{\varphi}$ by Lemma
\ref{le:1.4} (3) and (4).  Then, if we apply Theorem \ref{th:2.4} to both the
star operations $(\ast_{\varphi})^{\varphi}$ and $\ast$ on $R$, we
have the following chain of canonical isomorphisms:
$$
\Cl^{(\ast_{\varphi})^{\varphi}}(R)\, \cong \,
\Cl^{((\ast_{\varphi})^{\varphi})_{\varphi}}(D)\
=\ \Cl^{\ast_{\varphi}}(D)\, \cong\, \Cl^\ast(R)\,.
$$
Since these isomorphisms are canonical and $\Cl^\ast(R)$ is a subgroup
of $\Cl^{(\ast_{\varphi})^{\varphi}}(R)$ (Lemma \ref{le:2.2}), we easily conclude
that $\Cl^{(\ast_{\varphi})^{\varphi}}(R) = \Cl^\ast(R)$.
\end{proof}

\begin{reem}
\bf (1) \rm We  present  an example of a pullback diagram of type
$(\square^{+})$ in which $T$ is quasilocal and $\ast\lneq
(\ast_\varphi)^\varphi$ (with $\Cl^\ast\!(R)
=\Cl^{(\ast_{\varphi})^{\varphi}}\!(R)$ by  Corollary \ref{cor:2.7} (2)).
Let $D$ be an integral domain in which
each nonzero ideal is divisorial (e.g., a Dedekind domain) and let
$k$ be the quotient field of $D$. Set $T:=k[X^2,X^3]_Q$, where
       $Q:=X^2k[X]$, and $M:=QT$. Let $\varphi$ and $R$ be as in
$(\square^{+})$. Then
$((d_R)_\varphi)^\varphi=(d_D)^\varphi=(v_D)^\varphi=v_R$
\cite[Proposition 3.3 and Corollary 2.13]{fp}. Meanwhile, since
$T^{v_R}=(R:(R:T))=(R:M)\supseteq k[X]$ but $T\not\supseteq k[X]$,
$d_R\neq v_R=((d_R)_\varphi)^\varphi$.

\bf (2) \rm  We give an example to show that the quasilocal hypothesis is
essential in  Corollary \ref{cor:2.7} (2).  Let $D$ be an integral domain in which
each nonzero ideal is divisorial and let $k$ be the quotient field
of $D$. Let $B$ be the polynomial ring $k[\{X_i\}_{i=1}^\infty]$
and let $T$ be the subring of $B$ generated over $k$ by the
products $X_iX_j$ for all pairs $i,j \geq 1$. Then it is known that  $T$ is a Krull domain
\cite[Example 1.10]{fo}.
 Let $N:=(1+X_1, X_2,X_3,\cdots)B$ and let
$M:=N\cap T$. Since $k\subseteq T/M\subseteq B/N\cong k$,
$T/M\cong k$ and $T=k+M$. Let $\varphi$ and $R$ be as in
$(\square^{+})$. Then
$((d_R)_\varphi)^\varphi=(d_D)^\varphi=(v_D)^\varphi=v_R$
\cite[Proposition 3.3 and Corollary 2.13]{fp}. Let $Q:=X_1B\cap
T$ and note that $X_1B \ (\not\subseteq N)$ is a prime ideal of
height one in the Krull domain $B$.
       Since $B$ is integral over the integrally closed domain $T$,
 $Q$ is a prime ideal of
height one in $T$.
Note that $Q\not\subseteq M$,  because    $X_1^2\in Q\setminus N$.
Since $R =D+M$, $T=R_{D\setminus \{0\}}$, thus
$Q=\boldsymbol{q}T$, where $\boldsymbol{q}:=Q\cap R$ and $\boldsymbol{q} \not\subseteq M$.
 Since $Q$ is a prime ideal of
height one in the Krull domain, $Q$ is
a $t_T$--invertible $t_T$--ideal of $T$, thus $\boldsymbol{q}$ is a $t_R$--invertible
$t_R$--ideal of $R$ by \cite[ Lemma 3.1 and  Theorem 2.2 (6)]{ac}. Moreover, since $Q$
is not finitely generated as an ideal of $T$ \cite[Example 1.10]{fo}, $\boldsymbol{q}$ is not
finitely generated as an ideal of $R$ and hence it is not invertible. Therefore
$\Pic(R)=\Cl^{d_R}(R)\subsetneq \Cl^{t_R}(R)\subseteq \Cl^{v_R}(R)$,
thus $\Cl^{d_R}(R)\neq
\Cl^{v_R}(R)=\Cl^{((d_R)_\varphi)^\varphi}\!(R)$.

       This example also shows that the quasilocal hypothesis is essential in Corollary \ref{cor:2.7} (1):
Choose $D$ to be a PID. Then $\Cl^{d_D}(D)=\Pic(D)=0$, but since
$\Cl^{d_R}(R)\subsetneq \Cl^{v_R}(R)=\Cl^{(d_D)^\varphi}(R)$,
 we have   $\Cl^{(d_D)^\varphi}(R)\neq 0$.

\end{reem}

The next goal is to give a complete description of $\Cl^\ast(R)$ by
means of $\Cl^{\ast_{\varphi}}(D)$ and of an ``appropriate star class
group'' of $T$.  For this purpose, recall that, in \cite{fp}, we also
considered the problem of ``extending  a star operation'' defined on an
integral domain $R$ to  some overring   $T$ of $R$.

 We  need
the following notation. Let $\ast$ be a star operation on an integral
domain $R$ and let $T$ be
an overring of $R$ such that $(R:T) \neq 0$. Then, for each
$E \in \boldsymbol{F}(T)\ (\subseteq \boldsymbol{F}(R))$, we set
$$
E^{(\ast)_{_{\!T}}}:= E^\ast \cap (T:(T:E))= E^\ast \cap
E^{v_{T}}\,.
$$

  \begin{leem} \label{le:2.8}
Assume  that we are dealing with a pullback diagram of type
$(\square^{+})$.  Let $\iota:R \hookrightarrow T$ be the canonical
embedding and let $\ast$ be a star operation on $R$.
%

 \begin{enumerate}
\item[(1)] ${(\ast)_{_{\!T}}}$ is a star operation on $T$ with
${(\ast)_{_{\!T}}}= \ast_{\iota} \wedge v_{T}$.

\item[(2)] If $\ast$ is a star operation of finite type on $R$, then
${(\ast)_{_{\!T}}}$ coincides with ${\ast_{\iota}}$ (restricted to
the fractional ideals of $T$) and it is a star
operation of finite type on $T$.
\item[(3)] If $\ast_{1}, \ast_{2}$ are two star operations on $R$,
then
$$ \ast_{1} \leq \ast_{2} \;\;\; \Rightarrow \;\;\;
(\ast_{1})_{_{\!T}}\leq (\ast_{2})_{_{\!T}}\,.$$
\item[(4)] $(\ast_{_{\!f}})_{_{\!T}} \leq
\left((\ast)_{_{\!T}}\right)_{_{\!f}}$.
\item[(5)] $
(\widetilde{\ast})_{_{\!T}} = \widetilde{\ (\ast)_{_{\!T}}\ }$.
\end{enumerate}
\end{leem}
\begin{proof} (1) follows from  \cite[Example 1.2 and 1.5(a)]{fp}  and the observation
that $T^{(\ast)_{_{\!T}}} = T^\ast \cap T^{v_{T}} = T^\ast \cap T =T$.

For (2), we need the following:

\bf Claim 1.  \sl $T$ is a $t_{R}$--ideal of $R$.  \rm

Choose a nonzero $r \in M$, then obviously $rT$ is an integral $t_{T}$--ideal
of $T$ and $rT\subseteq M \subset R$.
Since $T$ is $R$--flat, $rT = rT \cap R$ is a $t_{R}$--ideal of $R$
by \cite[Proposition 0.7 (a)]{fg}.  Therefore, $T = r^{-1}\cdot rT$ is a
$t_{R}$--ideal of $R$.

By using Claim 1, we can complete the proof of (2).  As a
matter of fact, if $\ast$ is a star operation of finite type on $R$,
then $\ast \leq t_{R}$, thus the map $ E \mapsto E^{\ast_{\iota}}:=
E^\ast$, for each $E \in \boldsymbol{F}(T)\ (\subseteq
\boldsymbol{F}(R))$, defines a star operation on $T$ (since $T \subseteq
T^\ast \subseteq T^{t_{R}} =T$). In particular, $\ast_{\iota} \leq v_{T}$,
and so ${(\ast)_{_{\!T}}}=\ast_{\iota}$ (being $\ast_{\iota}$ restricted to
the fractional ideals of $T$).  Finally, it is straightforward
that if $\ast$ is a star operation of finite type on $R$, then
$\ast_{\iota}\ (= {(\ast)_{_{\!T}}})$ is of finite type on $T$ (cf.  also
for instance \cite[Example 1.2 (b)]{fp}).

(3) is a straightforward consequence of the definition.

(4) follows from (3) and (2) since $(\ast_{_{\!f}})_{_{\!T}}$ is a star
operation of finite type on $T$.

(5) Note that $(\widetilde{\ast})_{_{\!T}}$ is a star operation of
finite type and $(\widetilde{\ast})_{_{\!T}}= (\widetilde{\ast})_{\iota}$
(by (2)).  Moreover, $(\widetilde{\ast})_{\iota}$ is stable,
since $\widetilde{\ast}$ is stable.  Therefore
$(\widetilde{\ast})_{_{\!T}}= \widetilde{\ (\widetilde{\ast})_{_{\!T}}\
}$, and hence we conclude by (3) that $
(\widetilde{\ast})_{_{\!T}} \leq \widetilde{\ (\ast)_{_{\!T}}\ }$.

  \bf Claim 2.  \sl For each
star operation $\star$ on $R$, $M = M^{\star_{_{\!f}}} = M^\star$.  \rm

It follows from the fact that $M =(R:T)$ is a divisorial ideal of $R$.

\bf Claim 3.  \sl
 $\Max^{({\ast}_{_{\!f}})_{_{\!T}}}(T)$  coincides with the set of
maximal elements of $\{PT \mid P \in \Spec^{{\ast}_{_{\!f}}}(R)\,, \ PT
\neq T \}$.
 \rm

Since $T$ is $R$-flat  \cite[Lemma 0.3]{fg},
each ideal of $T$ is extended from $R$. In particular, each prime ideal $Q$ of
$T$ is equal to $(Q\cap R)T$.
Note that $\Max^{({\ast}_{_{\!f}})_{_{\!T}}}(T)
\subseteq \{PT \mid P \in \Spec^{{\ast}_{_{\!f}}}(R)\,, \ PT \neq T \}$.
Indeed, let $Q \in \Max^{({\ast}_{_{\!f}})_{_{\!T}}}(T)$
and let $P:=Q\cap R$. Then
$P \subseteq P^{{\ast}_{_{\!f}}}\subseteq Q^{{\ast}_{_{\!f}}}
=Q^{({\ast}_{_{\!f}})_{_{\!T}}} =Q $,
hence $P \subseteq P^{{\ast}_{_{\!f}}} \subseteq Q\cap R = P$.

Now let $PT$ be a maximal element in the set $ \{PT \mid P \in
\Spec^{{\ast}_{_{\!f}}}(R)\,, \ PT \neq T \}$.
Suppose $(PT)^{({\ast}_{_{\!f}})_{_{\!T}}}=T$.
Then $1\in (PT)^{({\ast}_{_{\!f}})_{_{\!T}}}=(PT)^{{\ast}_{_{\!f}}}$,
i.e., $1\in F^\ast$ for some $F\in \boldsymbol{f}(R)$ such that $F\subseteq PT$.
Let $m\in M\setminus \{0\}$. Then $m\in mF^\ast=(mF)^\ast\subseteq (mPT)^{{\ast}_{_{\!f}}}
\subseteq (PR)^{{\ast}_{_{\!f}}}=P^{{\ast}_{_{\!f}}}=P$.
Thus we have $M\subseteq P$.
Since $PT\neq T$, $M\not\subset P$, and hence $M=P$.
Then $T=(PT)^{({\ast}_{_{\!f}})_{_{\!T}}}=M^{({\ast}_{_{\!f}})_{_{\!T}}}
=M^{{\ast}_{_{\!f}}}=M$ (Claim 2), a contradiction.
Therefore, $(PT)^{({\ast}_{_{\!f}})_{_{\!T}}}\neq T$.

Let $Q' \in \Max^{({\ast}_{_{\!f}})_{_{\!T}}}(T)$
such that $(PT)^{({\ast}_{_{\!f}})_{_{\!T}}}\subseteq Q'$.
Then by the above argument, $Q'\cap R\in \Spec^{{\ast}_{_{\!f}}}(R)$.
Since $PT\subseteq Q'=(Q'\cap R)T$, $PT=Q'$ by the maximality of $PT$.
Thus we have $PT\subseteq (PT)^{({\ast}_{_{\!f}})_{_{\!T}}}\subseteq
Q'=PT$ and so $PT\in \Max^{({\ast}_{_{\!f}})_{_{\!T}}}(T)$.

\bf Claim 4. \sl    $\Max^{({\ast}_{_{\!f}})_{_{\!T}}}(T)=\Max^{((\ast)_{_{\!T}})_{_{\!f}}}(T)$.
\rm

Let $Q\in \Max^{((\ast)_{_{\!T}})_{_{\!f}}}(T)$ and let $P:=Q\cap R$.
  Then $P\subseteq
P^{{\ast}_{_{\!f}}}\subseteq Q^{{\ast}_{_{\!f}}}
=Q^{({\ast}_{_{\!f}})_{_{\!T}}}\subseteq
Q^{((\ast)_{_{\!T}})_{_{\!f}}}=Q$ (by (4)), and hence $P\subseteq
P^{{\ast}_{_{\!f}}}\subseteq Q\cap R=P$.   Thus we have
$\Max^{((\ast)_{_{\!T}})_{_{\!f}}}(T)\subseteq \{PT \mid P \in
\Spec^{{\ast}_{_{\!f}}}(R)\,, \ PT \neq T \}$.

Now let $PT$ be a maximal element in the set
$\{PT \mid P \in \Spec^{{\ast}_{_{\!f}}}(R)\,, \ PT \neq T \}$.
Suppose $(PT)^{((\ast)_{_{\!T}})_{_{\!f}}}=T$.
Then $1\in (PT)^{((\ast)_{_{\!T}})_{_{\!f}}}$, i.e.,
$1\in G^{(\ast)_{_{\!T}}}$ for some $G\in \boldsymbol{f}(T)$ such that
$G\subseteq PT$.     We may assume that $G=JT$ for some $J\in
\boldsymbol{f}(R)$ such that $J\subseteq P$.    Let $m\in M\setminus
\{0\}$.  Then $m\in
mG^{(\ast)_{_{\!T}}}=(mG)^{(\ast)_{_{\!T}}}=(mJT)^{(\ast)_{_{\!T}}}
\subseteq (mJT)^{\ast_\iota}=(mJT)^\ast\subseteq (JR)^\ast=J^\ast
\subseteq P^{{\ast}_{_{\!f}}}=P$.  Thus we have $M\subseteq P$.  Since
$PT\neq T$, $M\not\subset P$, and hence $M=P$.  Then
$T=(PT)^{((\ast)_{_{\!T}})_{_{\!f}}}=M^{((\ast)_{_{\!T}})_{_{\!f}}}
\subseteq M^{(\ast)_{_{\!T}}}\subseteq M^{\ast_\iota}=M^\ast=M$
(Claim 2),   a
contradiction.  Therefore, $(PT)^{((\ast)_{_{\!T}})_{_{\!f}}}\neq T$.

Let $Q'\in \Max^{((\ast)_{_{\!T}})_{_{\!f}}}(T)$
such that $(PT)^{((\ast)_{_{\!T}})_{_{\!f}}}\subseteq Q'$.
Then since $PT\subseteq Q'=(Q'\cap R)T$ and since we have already
proved that $Q'\cap R\in \Spec^{{\ast}_{_{\!f}}}(R)$, we conclude that
$PT=Q'$ by the maximality of $PT$.  Thus $PT\subseteq
(PT)^{((\ast)_{_{\!T}})_{_{\!f}}}\subseteq Q'=PT$ and so
$PT \in  \Max^{((\ast)_{_{\!T}})_{_{\!f}}}(T)$.

\bf Claim 5.  (a)  \sl For each prime ideal $P$ of $R$ such that $P
\not\supseteq M$, $R_{P} =TR_P=T_{PT}$; \bf (b) \sl for each prime
ideal $P$ of $R$ such that $P \supseteq M$, $R_P\subseteq R_M=T_M$, and
moreover, $TR_P=T_M$.  \rm

The statement (a) and the first part of (b) are well known \cite[Theorem
1.4 and its proof]{f}.  Since $TR_P \subseteq T_M$ for each $P \in\Spec(R)$
with $P \supseteq M$, to prove the equality, it suffices to show that if a
prime ideal $Q'$ of $T$ is such that $Q'\cap R\subseteq P$, then $Q'$
is contained in $M$.  Suppose not, i.e., $Q'\not\subseteq M$, then
$Q'\cap R \not\subseteq M$.
Choose $a\in (Q'\cap R)\setminus M$.
Then $M+aT=T$, so $1=m+at$ for some $m\in M$, $t\in T$.  Then $1-m=at\in
aT\cap R\subseteq Q'\cap R\subseteq P$.  Since $m\in M\subseteq P$,
$1\in P$, a contradiction.

 \bf Claim 6.  \sl
$\Max^{({\ast}_{_{\!f}})_{_{\!T}}}(T)
=
\{PT \mid P \in \Max^{{\ast}_{_{\!f}}}(R)\,, \ P \not\supseteq M\} \cup
\{M\}\,.  $ \rm

Note that, the condition $PT \neq T $ (or, equivalently, $PT \in \Spec(T)$)
implies that $P \not\supset M$, since $M$ is a maximal ideal in $T$.
Moreover, by Claim 2, $M$ belongs to $\Spec^{{\ast}_{_{\!f}}}(R)$, thus
$MT=M$ belongs, in any case, to $\Max^{({\ast}_{_{\!f}})_{_{\!T}}}(T)$ by
Claim 3.

Recall that, by the properties of the prime ideals in a pullback of type
$(\square^{+})$,  it follows that the canonical map $\Spec(T)
\rightarrow \Spec(R)$  is an order preserving embedding,
and if $Q \in \Spec(T)$ and  $Q\cap R \subseteq P$ for some $P \in \Spec(R)$ with $P
\supseteq M$, then  $Q  \subseteq M$  (see also  the proof of Claim 5).
By the previous ordering properties and Claim 3, we easily conclude that $
\{PT \mid P \in \Max^{{\ast}_{_{\!f}}}(R)\,, \ P \not\supseteq M\} \cup
\{M\} = \Max^{({\ast}_{_{\!f}})_{_{\!T}}}(T) $.

\bf Claim 7. \sl $\widetilde{\ (\ast_{_{\!f}})_{_{\!T}}\ }=(\widetilde{\ast})_{_{\!T}}$.
\rm

Note that, by  Claim 4,  $\widetilde{\ (\ast_{_{\!f}})_{_{\!T}}\ }
=\widetilde{\,  \left((\ast)_{_{\!T}}\right)_{_{\!f}}\, }
=\widetilde{\  (\ast)_{_{\!T}}\ }$.  Now we want to show that
$\widetilde{\ (\ast_{_{\!f}})_{_{\!T}}\ } =(\widetilde{\ast})_{_{\!T}}$.

Set $\mathcal{P}_{1}^{{\ast}_{_{\!f}}} :=\{P \in
\Spec^{{\ast}_{_{\!f}}}(R) \mid P \not\supseteq M\}$ and
$\mathcal{P}_{2}^{{\ast}_{_{\!f}}} :=\{P \in
\Spec^{{\ast}_{_{\!f}}}(R) \mid P \supseteq M\}$.
If we let $\mathcal{P}_{0}^{{\ast}_{_{\!f}}}$ be the set of maximal elements in the set
$\mathcal{P}_{1}^{{\ast}_{_{\!f}}}$, then $\{PT\mid P\in
\mathcal{P}_{0}^{{\ast}_{_{\!f}}}\} =
\{Q \in
\Max^{({\ast}_{_{\!f}})_{_{\!T}}}(T) \mid Q \neq M\}$ by Claim 6.

Let $E\in \boldsymbol{F}(T)$, then by using Claim 5 and 6, we have
$$
\begin{array}{rl}
    E^{(\widetilde{\ast})_{_{\!T}}} =& E^{(\widetilde{\ast})_{\iota}}  =
    E^{\widetilde{\ast}}=(ET)^{\widetilde{\ast}}
    = \cap \{ETR_{P} \mid P \in \Spec^{{\ast}_{_{\!f}}}(R) \}\\
    =& \left(\cap \{ETR_{P} \mid P \in \mathcal{P}_{1}^{{\ast}_{_{\!f}}} \}\right)
    \cap \left(\cap \{ETR_{P} \mid P \in
    \mathcal{P}_{2}^{{\ast}_{_{\!f}}} \}\right)\\
    =& \left( \cap \{ETR_P \mid P \in \mathcal{P}_{0}^{{\ast}_{_{\!f}}} \}\right)
        \cap ET_M \\
    =& \cap \{ET_{PT} \mid P \in  \Max^{{\ast}_{_{\!f}}}(R)\,,\; P \not\supseteq M \}\cap ET_{M}\\
     =&  \cap \{ET_{Q} \mid Q \in
       \Max^{(\ast_{_{\!f}})_{_{\!T}}}(T) \} \\
        =& E^{\widetilde{\ (\ast_{_{\!f}})_{_{\!T}}\ }}\,.  \\
\end{array}
$$

\vskip -18pt  \end{proof}

     \begin{reem}  \label{rk:2.12}
\bf (1) \rm  We were not able to prove or disprove the equality in the statement (4)
of Lemma \ref{le:2.8}. However
$(\ast_{_{\!f}})_{_{\!T}}=\left((\ast)_{_{\!T}}\right)_{_{\!f}}$
for the case $\ast=v_R$, which is the most important star operation
of nonfinite type. More precisely, \sl in the situation of Lemma
\ref{le:2.8}, we have
$$ (t_R)_{_{\!T}}= ((v_R)_{_{\!f}})_{_{\!T}}=\left((v_R)_{_{\!T}}\right)_{_{\!f}}\,.$$ \rm

Since $(t_R)_{_{\!T}}\leq \left((v_R)_{_{\!T}}\right)_{_{\!f}}$
and both terms are star operations of finite type (Lemma \ref{le:2.8}
(2)),
it suffices to show that $H^{(t_R)_{_{\!T}}}\supseteq H^{(v_R)_{_{\!T}}}$
for all nonzero finitely generated integral ideals $H$ of $T$.
Let $H$ be a nonzero finitely generated integral ideal of $T$.
Then $H=IT$ for some finitely generated ideal $I$ of $R$.

 If  $IT_M$ is not principal, then $I^{v_R}=I^{v_R}T$ by \cite[Proposition
2.7(1b)]{gh}.
Therefore, $H^{(v_R)_{_{\!T}}}\subseteq H^{(v_R)_\iota}=(IT)^{v_R}=
(I^{v_R}T)^{v_R}=I^{v_R}=I^{t_R}\subseteq
H^{t_R}=H^{(t_R)_\iota}=H^{(t_R)_{_{\!T}}}$.

 Now  assume that $IT_M$ is principal.
Then $H^{v_T}\subseteq (HT_M)^{v_{_{T_M}}}=(IT_M)^{v_{_{T_M}}}=IT_M$.
Let $R(M)$ be the CPI--extension of $R$ with respect to $M$, i.e.,
$R(M)$ is defined by the following pullback diagram \cite{bs}:
\[
\begin{array}{ccc}
R(M):={\varphi}^{-1}(D)& \longrightarrow & D \\
 \Big\downarrow & & \Big\downarrow
\\ T_M & \stackrel {\varphi}\longrightarrow & k=T_M/MT_M.
\end{array}
\]
Then by \cite[Lemma 1.3]{fg}, $R=R(M)\cap T$.
Note first that $TR(M)=T_M$,
because $TR(M)=\cap\, \{TR(M)_{\bar N}\mid {\bar N}\in \Max \left(R(M)\right)\}
=\cap\, \{TR_N \mid N\in \Max (R) \mbox { such that } \linebreak
 N\supseteq M\}=T_M$ by Claim 5 (b) in the proof of Lemma \ref{le:2.8}.
Now by \cite[Theorem 2(4)]{dda}, $H^{(t_R)_{_{\!T}}}=
H^{t_R}\supseteq \left( HR(M) \right)^{t_{R(M)}}\cap (HT)^{t_T}
=\left(ITR(M)\right)^{t_{R(M)}}\cap H^{v_T}
=\left(IT_M\right)^{t_{R(M)}}\cap H^{v_T}\supseteq IT_M\cap H^{v_T}=H^{v_T}
\supseteq H^{(v_R)_{_{\!T}}}$.

\bf (2) \rm   As another special case,  we have the following positive result.
\sl Consider a pullback diagram of type $(\square^+)$,
let $\star'$ be a star operation on $D$ and  $\star''$   a star operation
on $T$.  Set \small{$\diamond$}\normalsize$:={\star'}^{\varphi} \wedge
{\star''}^{{\iota}}$.
We know that \small{$\diamond$} \normalsize  is a star
operation on $R$ \cite[Corollary 2.5]{fp}.
If  $(({\star'}^{\varphi})_{_{\!T}})_{_{\!f}}  =
(({\star'}^{\varphi})_{_{\!f}})_{_{\!T}}$  (e.g. this hypothesis
is satisfied in each one of the following cases: \rm\bf (a)  \sl $\star'=
v_{D}$,\  \rm\bf  (b) \sl   $({\star'_{_{\!f}}})^{\varphi}$
is a star operation of finite type  on $R$,\  \rm\bf  (c) \sl  $T$ is
a Pr\"{u}fer domain),
then $
(\mbox{\small {$\diamond$}}_{_{\!f}})_{_{\!T}} =
((\mbox{{\small{$\diamond$}}})_{_{\!T}})_{_{\!f}} $. \rm

\smallskip

\bf Claim 1. \sl  If $\ast_{1}$ and  $\ast_{2}$ are two  semistar
operations on an integral domain $R$, then $(\ast_{1}\wedge
\ast_{2})_{_{\!f}} = (\ast_{1})_{_{\!f}} \wedge
(\ast_{2})_{_{\!f}}$\,. \rm

 This is an easy consequence of the
fact that ``$\bigcup_{\alpha}$
distributes over $\cap$''.
\smallskip

\bf Claim 2. \sl Let  $\iota: R \hookrightarrow T$ be  an embedding of
an integral domain $R$ in one of its overrings $T$ and let $\star$ be a  semistar
operation on $T$.  Then, in $R$,  $(\star^\iota)_{_{\!f}}  =
(\star_{_{\!f}})^\iota$, and in $T$,   $\star = (\star^{\iota})_\iota$
(Example 1.1(e3))\,. \rm

Let  $E \in \boldsymbol{\overline{F}}(R)$ and let $G \in
    \boldsymbol{f}(T)$ be contained in $ET$. Then  $G :=(x_{1}t_{1},
    x_{2}t_{2}, \ldots, \linebreak x_{n}t_{n})T$ for some $n\geq 1$, $\{x_{1},
    x_{2}, \ldots, x_{n} \} \subseteq E$, and  $\{t_{1},
    t_{2}, \ldots, t_{n} \} \subseteq T$.  Thus $G \subseteq HT$,
    where $H := (x_{1},
    x_{2}, \ldots, x_{n})R\in \boldsymbol{f}(R)$ (and $H \subseteq E$). Therefore
$$
\begin{array}{rl}
    E^{(\star^\iota)_{_{\!f}}} =& \cup \{F^{\star^\iota} \mid F \in
    \boldsymbol{f}(R)\,,\; F \subseteq E \} \\
    =& \cup \{(FT)^{\star} \mid F \in
    \boldsymbol{f}(R)\,,\; F \subseteq E \}\\
    = &  \cup \{G^{\star} \mid G \in
    \boldsymbol{f}(T)\,,\; G \subseteq ET \}\\
    =& (ET)^{\star_{_{\!f}}} = E^{(\star_{_{\!f}})^\iota} \,.
    \end{array}
    $$

\bf Claim 3. \sl Let  $\iota: R \hookrightarrow T$ be  an embedding of
an integral domain $R$ in one of its overrings $T$ and let  $\ast_{1}$
and  $\ast_{2}$ be two  semistar
operations  on $R$.  Then $(\ast_{1}\wedge
\ast_{2})_{\iota} = (\ast_{1})_{\iota} \wedge
(\ast_{2})_{\iota}$\,. \rm

This is an obvious consequence of the definitions.

\smallskip

\bf Claim 4. \sl Let  $\iota: R \hookrightarrow T$ be  an embedding of
an integral domain $R$ in one of its overrings $T$ and let $\ast$ be
a semistar operation on $R$. Then $(\ast_{_{\!f}})_{\iota}$ is a
semistar operation of finite type on $T$. \rm

For each  $E\in \boldsymbol{\overline{F}}(T)$,  we have
$$
\begin{array}{rl}
    E^{(\ast_{_{\!f}})_{\iota}} = &  E^{\ast_{_{\!f}}} = \cup \{F^{\ast} \mid F \in
    \boldsymbol{f}(R)\,,\; F \subseteq E \} \\
    =& \cup \{\ \cup  \{F^{\ast} \mid F \in
    \boldsymbol{f}(R)\,,\; F \subseteq G \}\  \mid G \in
    \boldsymbol{f}(T)\,,\; G \subseteq E \}\\
    =& \cup \{G^{\ast_{_{\!f}}} \mid G \in
    \boldsymbol{f}(T)\,,\; G \subseteq E \} \\
     =& \cup \{G^{(\ast_{_{\!f}})_{\iota}} \mid G \in
    \boldsymbol{f}(T)\,,\; G \subseteq E \} \\
    =&  E^{((\ast_{_{\!f}})_{\iota})_{_{\!f}}}\,.
    \end{array}
    $$

\bf Claim 5. \sl  In a pullback diagram of type $(\square)$,
let $\star$ be a star operation on $D$. Then
    $(\star^\varphi)_{\iota} = (v_{R})_{\iota}$\  (when restricted to
    $\boldsymbol{F}(T)$),
    and hence $(\star^\varphi)_{_{\!T}} = (v_{R})_{_{\!T}}$.
    Moreover, in a pullback diagram of type  $(\square^{+})$,
    $((\star^\varphi)_{_{\!T}})_{_{\!f}} = (t_{R})_{_{\!T}}$ by (1).  \rm

 Let $I$ be a nonzero integral ideal of $T$.  Note that
$$
x \in (R : I) \; \Rightarrow \; xIT =xI \subseteq R \; \Rightarrow \;
xI \subseteq (R:T) = M \; (\, \Leftrightarrow \;  I \subseteq
x^{-1}M\ ).
$$
Therefore we have
$$
\begin{array}{rl}
I^{(\star^\varphi)_{\iota}} = & I^{\star^{\varphi}} =
 \cap \left\{
  x^{-1}\varphi^{-1}\left(\left(\frac{xI+M}{M}\right)^{\star}\right) \mid
  \, x \in (R:I)\,, \; x \neq 0
  \right\} \\
  = &  \cap \left\{
  x^{-1}M \mid \, x \in (R:I)\,, \; x \neq 0
  \right\}  = I^{v_{R}} = I^{(v_{R})_{\iota}} \,.
  \end{array}
  $$
  Note that $T^{(\star^\varphi)_{\iota}} =  T^{(v_{R})_{\iota}}=
  T^{v_{R}}$, thus ${(\star^\varphi)_{\iota}}$ (when restricted to
    $\boldsymbol{F}(T)$) is a star operation
  on $T$ if and only if $T = T^{v_{R} }$.

\smallskip

Now we use the previous claims to prove the statement.
By applying Claim 2, 3, and 5, we have
$$
\begin{array}{rl}
   \mbox{\small{($\diamond$)}}_{_{\!T}} =&
    \mbox{\small{$\diamond$}}_{\iota} \wedge v_{T}=
    ({\star'}^{\varphi} \wedge
{\star''}^{{\iota}})_{\iota} \wedge v_{T} \\
=& ({\star'}^{\varphi})_{\iota} \wedge
({\star''}^{{\iota}})_{\iota} \wedge v_{T} =
({\star'}^{\varphi})_{\iota} \wedge
{\star''} \wedge v_{T}\\
=& ({\star'}^{\varphi})_{\iota} \wedge
{\star''} = (v_{R})_{\iota} \wedge
{\star''}  \;\; \; \mbox{ or equivalently, }\\
=&   ({\star'}^{\varphi})_{_{\!T}} \wedge
{\star''} = (v_{R})_{_{\!T}}   \wedge
{\star''} \,.
\end{array}
$$
Therefore, by Claim 1 and (1), we have
$$
 (\mbox{\small{($\diamond$)}}_{_{\!T}})_{_{\!f}}
 = (({\star'}^{\varphi})_{\iota})_{_{\!f}} \wedge
{\star''_{_{\!f}}} =
 (({\star'}^{\varphi})_{_{\!T}})_{_{\!f}} \wedge
{\star''_{_{\!f}}} =
((v_{R})_{_{\!T}})_{_{\!f}}  \wedge
{\star''_{_{\!f}}}  =(t_{R})_{_{\!T}}  \wedge
{\star''_{_{\!f}}} \,.
$$
On the other hand,
by Lemma \ref{le:2.8}(1), Claim 1, 2  and 3, we have
$$
\begin{array}{rl}
   (\mbox{\small{$\diamond$}}_{_{\!f}})_{_{\!T}} =&
   (\mbox{\small{$\diamond$}}_{_{\!f}})_{\iota} =
   (({\star'}^{\varphi})_{_{\!f}})_{\iota} \wedge
(({\star''}^{{\iota}})_{_{\!f}})_{\iota}  =
   (({\star'}^{\varphi})_{_{\!f}})_{\iota} \wedge
((\star''_{_{\!f}})^{{\iota}})_{\iota}  \\
=&
(({\star'}^{\varphi})_{_{\!f}})_{\iota}  \wedge
\star''_{_{\!f}}  = (({\star'}^{\varphi})_{_{\!f}})_{_{\!T}}  \wedge
\star''_{_{\!f}}\,.
\end{array}
$$
It is obvious now that, if  $(({\star'}^{\varphi})_{_{\!T}})_{_{\!f}}  =
(({\star'}^{\varphi})_{_{\!f}})_{_{\!T}}$, then  $
(\mbox{\small {$\diamond$}}_{_{\!f}})_{_{\!T}} =
((\mbox{{\small{$\diamond$}}})_{_{\!T}})_{_{\!f}}$.

 Finally, we check the parenthetical statement.

Assume that $\star' =v_{D}$, then we know that $(v_{D})^\varphi =
v_{R}$ \cite[Corollary 2.13]{fp}.  Therefore $(({\star'}^{\varphi})_{_{\!f}})_{_{\!T}}=
(t_{R})_{_{\!T}}$ and so $(({\star'}^{\varphi})_{_{\!f}})_{_{\!T}}$
coincides with $(({\star'}^{\varphi})_{_{\!T}})_{_{\!f}}=
((v_{R})_{_{\!T}})_{_{\!f}}$ by (1).

Assume that $({\star'_{_{\!f}}})^{\varphi} $  is a star operation of
finite type. Note that, from the fact that
$ ({\star'_{_{\!f}}})^{\varphi} \leq {\star'}^{\varphi} $ and
from the
assumption, it follows that $ ({\star'_{_{\!f}}})^{\varphi} \leq
({\star'}^{\varphi})_{_{\!f}}$. Therefore, by \cite[Proposition 2.9, Theorem 2.12 and
Proposition 3.6(b)]{fp}, we have
$$
    ({\star'}^{\varphi})_{_{\!f}} \leq  ((({\star'}^{\varphi})_{_{\!f}})_{\varphi})^{\varphi}
=  ((({\star'}^{\varphi})_{\varphi})_{_{\!f}})^{\varphi}
= (\star'_{_{\!f}})^{\varphi}\,,
$$
thus $ ({\star'_{_{\!f}}})^{\varphi} =({\star'}^{\varphi})_{_{\!f}}$.
In this situation,
by Claim 5,
we have $ (t_{R})_{\iota}\leq (v_{R})_{\iota}=
(({\star'_{_{\!f}}})^{\varphi})_{\iota}=
(({\star'}^{\varphi})_{_{\!f}})_{\iota}\leq
(t_{R})_{\iota}$. Therefore, $
(t_{R})_{\iota} =
(({\star'_{_{\!f}}})^{\varphi})_{\iota}=(({\star'}^{\varphi})_{_{\!f}})_{\iota}=(v_{R})_{\iota}
=((v_{R})_{\iota})_{_{\!f}}$ and so, in particular,
$
(t_{R})_{_{\!T}} =
(({\star'_{_{\!f}}})^{\varphi})_{_{\!T}}=(({\star'}^{\varphi})_{_{\!f}})_{_{\!T}}=
(v_{R})_{_{\!T}}$\,. On the other hand, by Claim 5, we
know that $(({\star'}^{\varphi})_{_{\!T}})_{_{\!f}} =
((v_{R})_{_{\!T}})_{_{\!f}}= (t_{R})_{_{\!T}}$\,.

Assume that $T$ a Pr\"{u}fer  domain, then clearly $T$ has a unique
star operation of finite type, since $d_{T}= t_{T}$. In this situation,
obviously
$d_{T}=(({\star'}^{\varphi})_{_{\!f}})_{_{\!T}} =
(t_{R})_{_{\!T}}=t_{T}$, and
from Claim 5, we have
$(({\star'}^\varphi)_{_{\!T}})_{_{\!f}} = (t_{R})_{_{\!T}}$.

\bf (3) \rm  Under the assumptions of Lemma \ref{le:2.8}, as a consequence of Claim 3 and 6 in its proof,  we have
    that $\Max^{{(t_{R}){_{_{\!T}}}}}(T)$   coincides with the set of the maximal
         elements of $\{PT \in \Spec(T) \mid P\in \Spec^{t_{R}}(R)\} $ (which is
         equal to the set $\{PT \mid P\in \Max^{t_{R}}(R)\,, \ P\not\supseteq M\}
           \cup \{M\}$).

           We can give a little different proof of this result    under the additional assumption that
         the map $\tilde{\varphi}: \mathcal
U(T)\rightarrow k^\bullet/\mathcal{U}(D)$ is surjective.
Let $Q\in \Max^{{(t_{R}){_{_{\!T}}}}}(T)$ and let $P:=Q\cap R$.  Then $Q=PT$ and
         $Q=Q^{{(t_{R}){_{_{\!T}}}}}=Q^{t_R}$.  Therefore $P \subseteq
         P^{t_R}\subseteq Q^{t_R}\cap R= Q \cap R = P$ and so
         $\Max^{{(t_{R}){_{_{\!T}}}}}(T) \subseteq \{PT \in \Spec(T) \mid P\in
         \Spec^{t_{R}}(R)\}$.

         Conversely,   let $Q:=PT$ be a maximal element of the set $\{PT\in  \Spec(T)
         \mid P \in \Spec^{t_{R}}(R) \}$.    Assume that $P =M$, then    since $M =
         MT$ is a maximal ideal of $T$ and $ M = M^{t_{R}}$, $M$    is also
         a ${{(t_{R}){_{_{\!T}}}}}$--ideal of $T$, thus $M=MT\in
         \Max^{{{(t_{R}){_{_{\!T}}}}}}(T)$.  Assume that $P\neq M$.  Then
         $P\not\subseteq M$ by the maximality of $Q=PT$.  Now, if
         $\boldsymbol{S} :=\mathcal U(T)\cap R$, then by \cite[Theorem 2.2
         (5) and Lemma 3.1]{ac} we have $(PT)^{t_T}=(PR_{\boldsymbol
         S})^{t_T}=P^{t_R}R_{\boldsymbol S}=PR_{\boldsymbol S}=PT$.  Since
         $Q=PT\in \Spec^{t_{T}}(T)$,  $Q\in
         \Spec^{{{(t_{R}){_{_{\!T}}}}}}(T)$.

\end{reem}

\begin{leem} \label{le:2.9} Assume that we are dealing with a pullback diagram of
type $(\square^{+})$. Let $\ast$ be a
star operation of finite type on $R$ and let ${(\ast)_{_{\!T}}}$ be the star
operation on $T$ defined just before Lemma \ref{le:2.8}.
\begin{enumerate}
\item[(1)] If $H\in \Inv^\ast(R)$, then  $HT \in \Inv^{{(\ast)_{_{\!T}}}}(T)$.

\item[(2)] The canonical map $\boldsymbol{\beta}(\varphi, \ast)$ (or, simply,
$\boldsymbol{\beta}$): $\Inv^\ast(R)\rightarrow \Inv^{{(\ast)_{_{\!T}}}}(T)$,
$H\mapsto HT$, is a group-homomorphism.

\item[(3)] The map $\boldsymbol{\beta}$, defined in \rm (2), \it induces a
group-homomorphism $\overline{\boldsymbol{\beta}}(\varphi, \ast)$\ (or,
simply, $\overline{\boldsymbol{\beta}}$) : $\Cl^\ast(R) \rightarrow
\Cl^{{(\ast)}_{_{\!T}}}(T)$, $[H]\mapsto [HT]$.
\end{enumerate}
\end{leem}
\begin{proof} (1) Note that if $H$ is a $\ast$--invertible $\ast$--ideal of $R$
and $\ast = \ast_{_{\!f}}$, then $H$ is a $t_R$-invertible $t_R$--ideal
of $R$  (Lemma \ref{lemma:1.2} (2)).   Moreover,
$T$ is a
flat overring of $R$ \cite[Lemma 0.3]{fg},  and hence $HT$ is a $t_T$-invertible $t_T$-ideal of
$T$ \cite[Proposition 0.7 (b)]{fg}.  We know by Lemma \ref{le:2.8} (2)
that ${(\ast)_{_{\!T}}}$ is a star operation of finite-type on $T$, so
${(\ast)_{_{\!T}}}\leq t_T$, and hence $HT$ is a
${(\ast)_{_{\!T}}}$--ideal of $T$.  Now, we show that $HT$ is also
${(\ast)_{_{\!T}}}$--invertible:
$$
\begin{array}{rl}
    (HT(HT)^{-1})^{{(\ast)_{_{\!T}}}} = & (HT(HT)^{-1})^{\ast}\cap
    (HT(HT)^{-1})^{v_T}=(HT(HT)^{-1})^\ast \cap T \\
    \supseteq & (HH^{-1}T)^\ast \cap T=((HH^{-1})^\ast T)^\ast \cap T
    \\
    = & (RT)^\ast \cap T =T^\ast \cap T=T\,, \end{array}
    $$
    thus $1 \in (HT(HT)^{-1})^{{(\ast)_{_{\!T}}}}$ and so $T
    =(HT(HT)^{-1})^{(\ast){_{_{\!T}}}}$.

   (2) is an obvious consequence of (1) and (3) follows from (2).
    \end{proof}

\begin{thee}\label{thm1}
Assume that we are dealing with a pullback diagram of type
$(\square^{+})$.  Suppose that the map $\tilde{\varphi}: \mathcal
U(T)\rightarrow k^\bullet/\mathcal{U}(D)$ is surjective and that $\ast$ is
a star operation of finite type on $R$.
Then
$\overline{\boldsymbol{\beta}}:=\overline{\boldsymbol{\beta}}(\varphi, \ast) : \Cl^\ast(R) \rightarrow
\Cl^{{(\ast)}_{_{\!T}}}(T)$ is surjective.
\end{thee}
\begin{proof}
Let $J$ be an integral
$(\ast){_{_{\!T}}}$--invertible $(\ast){_{_{\!T}}}$--ideal of $T$.  Then
$J=(IT)^{(\ast)_T}=(IT)^{t_T}$ for some finitely generated integral ideal $I$ of $R$
(\cite[Proposition 3.1 and Proposition 3.2]{a} and \cite[Lemma 0.3]{fg}). 

\bf Claim 1. \sl Without loss of generality, we may assume that $I\not\subseteq M$ . \rm

Suppose that $II^{-1}\subseteq M$. Then
$$
\begin{array}{rl}
    \left(JJ^{-1}\right)^{(\ast)_{_{\!T}}}
    =&
    \left((IT)^{{(\ast)_{T}}}\left((IT)^{{(\ast)_{T}}}\right)^{-1}\right)^{(\ast)_{T}}\\
    =& \left((IT)(IT)^{-1}\right)^{(\ast)_{T}}  \\
     =&
     \left(II^{-1}T\right)^{(\ast)_{T}} \\
      \subseteq & (MT)^{(\ast)_{T}}=M^{(\ast)_T}=M \,,
    \end{array}
    $$
which contradicts that $J$ is $(\ast)_T$-invertible.
    Thus, $II^{-1}\not\subseteq M$ and so we can choose  $x\in I^{-1}$  such that $xI\not\subseteq M$.
Set $I':=xI$ and $J':=xJ$. Then $I'\not\subseteq M$ and $J'=(I'T)^{(\ast)_T}$.
Since the classes $[J]$ and $[J']$ in $\Cl^{({{\ast})}_{_{\!T}}}(T)$ are the same, we can replace $J$ by
$J'$ and $I$ by $I'$.

Set $\boldsymbol{S} :=\mathcal U(T)\cap R$ (as in Remark \ref{rk:2.12}) and
  $\boldsymbol{N}:=\{x\in R\mid \varphi(x)\in\mathcal{U}(D)\}$. Then
$T=R_{\boldsymbol{S}}$ and $\boldsymbol{S}\cdot\boldsymbol{N}
=R\setminus M$ \cite[Lemma 3.1]{ac}. Since we may assume that
$I\not\subseteq M$, by \cite[Theorem 2.2 (2)]{ac} we have $I^{t_R} =
((S_{1})(N_{1}))^{t_R}$ for some nonempty finite subsets $S_1$ of
$\boldsymbol{S}$ and $N_1$ of $\boldsymbol{N}$.  Again by
\cite[Theorem 2.2]{ac},
$J=(IT)^{t_T}=I^{t_R}T=((S_1)(N_1))^{t_R}T=((S_1)(N_1)T)^{t_T}=((N_1)T)^{t_T}=
(N_1)^{t_R}T$, and hence $JJ^{-1}=
((N_1)^{t_R}T)(((N_1)T)^{t_T})^{-1}=((N_1)^{t_R}T)((N_1)T)^{-1}=
(N_1)^{t_R}(N_1)^{-1}T$.

\bf Claim 2. \sl If $\ast=\widetilde{\ast}$, then $\overline{\boldsymbol{\beta}}$ is surjective. \rm

  Let $P'\in \Spec^\ast(R)$ such that $M\not\subseteq P'$.  Then there
exists a unique prime ideal $Q'$ of $T$ such that $Q'\cap R=P'$ and
$R_{P'}=T_{Q'}$ \cite[Theorem 1.4, point (c) of the proof]{f}.  Since
$T=(JJ^{-1})^{(\ast){_{_{\!T}}}}=(JJ^{-1})^{\ast_\iota}=\cap
\{JJ^{-1}R_P \mid P\in \Max^\ast (R)\}= \cap \{JJ^{-1}R_P \mid P\in
\Spec^\ast (R)\} \subseteq JJ^{-1}R_{P'}=JJ^{-1}T_{Q'}$,
$JJ^{-1}\not\subseteq Q'$, and hence
$(N_1)^{t_R}(N_1)^{-1}\not\subseteq P'$.

Now let $P''\in \Spec^\ast(R)$ such that $M\subseteq P''$.  Then $P'' \cap
\boldsymbol{N}=\emptyset$, because if $x\in P''\cap \boldsymbol{N}$, then $\varphi(x)\in
P''/M\in \Spec(D)$, which contradicts that $\varphi(x)\in \mathcal{U}(D)$.
Therefore $(N_1)^{t_R}(N_1)^{-1}\not\subseteq P''$.

Thus since $(N_1)^{t_R}(N_1)^{-1}\not\subseteq P$ for all $P\in \Spec^\ast (R)$,
$((N_1)^{t_R}(N_1)^{-1})^\ast = R$, i.e., $(N_1)^{t_R}$ is a
$\ast$--invertible $\ast$--ideal of $R$.  Therefore, passing to the
classes,
$[J]=[(N_1)^{t_R}T]=\overline{\boldsymbol{\beta}}([(N_1)^{t_R}])$.

 \bf Claim 3. \sl $\Cl^{{(\ast){_{_{\!T}}}}}(T)=\Cl^{{(\widetilde{\ast}){_{_{\!T}}}}}(T)$ \rm
(it does hold without the condition $\tilde{\varphi}: \mathcal
U(T)\rightarrow k^\bullet/\mathcal{U}(D)$ is surjective).

By \cite[Theorem 2.18]{ack} and Lemma \ref{le:2.8} (5),
$\Cl^{{(\ast){_{_{\!T}}}}}(T)
 =\Cl^{\widetilde{\ ({\ast}){_{_{\!T}}}\ }}(T)= \Cl^{{(\widetilde{\ast}){_{_{\!T}}}}}(T)$.

Finally, since $\Cl^{\ast}(R)=\Cl^{\widetilde{\ast}}(R)$ by \cite[Theorem 2.18]{ack},
$ \overline{\boldsymbol{\beta}}(\varphi, \ast) =
\overline{\boldsymbol{\beta}}({\varphi}, \widetilde{\ast})$ and hence the conclusion follows.
\end{proof}

From Claim 3 in the proof of Theorem ~\ref{thm1} we deduce
immediately:

\begin{coor}\label{cor2bis}
    Assume that we are dealing with a pullback diagram of type
$(\square^{+})$. Then
$\Cl^{{(t_{R}){_{_{\!T}}}}}(T)$ $=\Cl^{{(w_{R}){_{_{\!T}}}}}(T)$.
\end{coor}

    In order to give a description of $\Cl^\ast(R)$ by means of
    $\Cl^{\ast_{\varphi}}(D)$ and $\Cl^{{(\ast)}_{_{\!T}}}(T)$,
     we need the following result from  \cite{fg}:

    \begin{leem} \label{le:gamma} \rm (\cite[Lemma 2.2 and the subsequent
    considerations]{fg}) \it Assume that we are dealing with a pullback
    diagram of type $(\square^{+})$.
\begin{enumerate}
    \item[(1)] For each $H \in \Inv^t(R)$ there exist a nonzero element $z$
    in the quotient field of $R$ and $H' \in \Inv^t(R)$, with $H'
    \not\subseteq M,\ H' \subseteq R$,\ and $H = zH'$.

    \item[(2)] The map $\overline{\boldsymbol{\gamma}}: \Cl^{t_{R}}(R)
    \rightarrow \Cl^{t_{D}}(D)$,\ $[H] \mapsto [(\varphi(H'))^{v_{D}}]$, is a
    well-defined group-homomorphism (where $H'$ is chosen as in \rm (1)\it\
    ).
    \end{enumerate}
    \end{leem}

    \begin{coor} \label{cor:gamma}  Assume that we are dealing with a pullback
    diagram of type $(\square^{+})$.  Let $\overline{\boldsymbol{\gamma}}:
    \Cl^{t_{R}}(R) \rightarrow \Cl^{t_{D}}(D)$ be as in Lemma \ref{le:gamma}
    and let $\ast$ be a star operation of finite type on $R$.  Then, by
    restriction to $ \Cl^{\ast}(R) \ (\subseteq \Cl^{t_{R}}(R))$,
    $\overline{\boldsymbol{\gamma}}$ defines a group-homomorphism
    $\overline{\boldsymbol{\gamma}}
    =:\overline{\boldsymbol{\gamma}}(\varphi, \ast): \Cl^{\ast}(R)
    \rightarrow \Cl^{\ast_{\varphi}}(D)$.
    \end{coor}

    \begin{proof} We want to show that
    $\overline{\boldsymbol{\gamma}}(\Cl^\ast(R))\subseteq
    \Cl^{\ast_{\varphi}}(D)\subseteq \Cl^{t_D}(D)$.  First, recalling that
    $\ast_{\varphi} \leq (t_R)_{\varphi}= t_D$  \cite[Proposition 3.7]{fp},
    we have $\Cl^{\ast_{\varphi}}(D)\subseteq \Cl^{t_D}(D)$.  Now let $H$ be a
    $\ast$--invertible $\ast$--ideal of $R$ such that $H\subseteq R$ and
    $H\not\subseteq M$.  Choose $r\in H\smallsetminus M$.  Then
    $rH^{-1}\subseteq R$ and $rH^{-1}\not\subseteq M$.  By using the fact
    that $\varphi(r)D$ is a divisorial ideal of $D$ and \cite[Proposition
    2.7]{fp}, we have
    $$
    \begin{array}{rl} \varphi(r)D=& (\varphi(r)D)^{\ast_{\varphi}}=
    (\varphi(rR))^{\ast_{\varphi}}=\left(\varphi\left(r\left(HH^{-1}\right)^\ast\right)\right)^{\ast_{\varphi}}\\
    =&
    \left(\frac{r\left(HH^{-1}\right)^\ast+M}{M}\right)^{\ast_{\varphi}}=
    \frac{\left(r\left(HH^{-1}\right)^\ast+M\right)^\ast}{M}\\
    =& \frac{\left(rHH^{-1}+M\right)^\ast}{M}=\left(\frac{rHH^{-1}+M}{M}\right)^{\ast_{\varphi}}\\
    =& \left(\frac{H+M}{M}\frac{rH^{-1}+M}{M}\right)^{\ast_{\varphi}}
    =\left(\varphi(H)\varphi\left(rH^{-1}\right)\right)^{\ast_{\varphi}}\,.
    \end{array}
    $$
    Hence $\varphi(H)$ is $\ast_{\varphi}$--invertible, and so
    $(\varphi(H))^{v_D}$ is a $\ast_{\varphi}$--invertible $\ast_{\varphi}$--ideal of $D$
    (Lemma  \ref{lemma:inv} (5)).  Therefore
    $\overline{\boldsymbol{\gamma}}$ induces a homomorphism
    $\overline{\boldsymbol{\gamma}}(\varphi, \ast): \Cl^\ast(R)\rightarrow
    \Cl^{\ast_{\varphi}}(D)$.  \end{proof}

    \begin{thee}\label{th:split}   Assume that we are dealing with a pullback
    diagram of type $(\square^{+})$.  Suppose that the map $\tilde{\varphi}:
    \mathcal U(T)\rightarrow k^\bullet/\mathcal U(D)$ is surjective and that
    $\ast$ is a star operation of finite type on $R$.  Then
the sequence
$$ 0\rightarrow \Cl^{\ast_{\varphi}}(D) \stackrel
{\overline{\boldsymbol{\alpha}}} {\longrightarrow} \Cl^\ast(R) \stackrel
{\overline{\boldsymbol{\beta}}}
{\longrightarrow}\Cl^{{(\ast)}_{_{\!T}}}(T)\rightarrow 0$$
is split exact.
\end{thee}
\begin{proof} It is obvious that $\overline{\boldsymbol{\alpha}}$ is
injective, since $\boldsymbol{\alpha}$ is injective (Lemma
\ref{le:2.1}).
The surjectivity of $\overline{\boldsymbol{\beta}}$ follows from Theorem \ref{thm1}.
To see that $\Ima(\overline{\boldsymbol{\alpha}})
=\Ker(\overline{\boldsymbol{\beta}})$, let $[H] \in
\Ima(\overline{\boldsymbol{\alpha}})$.
We can assume that $H =\varphi^{-1}(J)$ for
some $J \in \Inv^{\ast_{\varphi}}(D)$ and so $M\subsetneq H \subseteq
T$.  Hence, in particular, $HT=T$, because $M$ is a maximal ideal of $T$,
and thus $\overline{\boldsymbol{\beta}}([H]) = [HT]=[T]$. \ Conversely, let $[H] \in
\Ker(\overline{\boldsymbol{\beta}})$.  Without loss of generality, we can
assume that $H\in \Inv^\ast(R)$ and $HT=T$.
Then by \cite[Proposition 1.1]{fg} and \cite[Proposition 3.1 (a)]{a},
$M\subsetneq H=H^{v_R}\subseteq T$.
Moreover, since $T$ is not a $\ast$--invertible ($\ast$--)ideal of $R$, $H^{v_R}\subsetneq T$.
By Lemma \ref{le:2.1},
$H={\varphi}^{-1}(J)$ for some $\ast_{\varphi}$--invertible
$\ast_{\varphi}$--ideal $J$ of $D$, hence $H\in
\Ima(\overline{\boldsymbol{\alpha}})$.  Thus the sequence is exact.

Lastly, by the definitions of  ${\overline{\boldsymbol{\alpha}}}
={\overline{\boldsymbol{\alpha}}}(\varphi, \ast)$ and
${\overline{\boldsymbol{\gamma}}}
={\overline{\boldsymbol{\gamma}}}(\varphi, \ast)$ (Lemma \ref{lemma2}
and Corollary \ref{cor:gamma}), we immediately obtain  that
${\overline{\boldsymbol{\gamma}}} \circ
{\overline{\boldsymbol{\alpha}}}: \Cl^{\ast_{\varphi}}(D) \rightarrow \Cl^{\ast}(R) \rightarrow
\Cl^{\ast_{\varphi}}(D)$ is such that $[J] \mapsto
{\overline{\boldsymbol{\gamma}}}\left(\left[\varphi^{-1}(J)\right]\right)
= \left[{\left(\varphi\left(\varphi^{-1}(J)\right)\right)}^{v_D}\right] =[J^{v_D}]=[J]$,
i.e., it is the identity map.  Therefore the above exact sequence
splits.
\end{proof}

\begin{coor} \label{cor:2.13}  Assume that we are dealing with a pullback
    diagram of type $(\square^{+})$ and that the map $\tilde{\varphi}:
    \mathcal U(T)\rightarrow k^\bullet/\mathcal U(D)$ is surjective.
Then the sequence
$$
0\longrightarrow \Cl^{t_{D}}(D) \stackrel {\overline{\boldsymbol{\alpha}}}
{\longrightarrow} \Cl^{t_{R}}(R) \stackrel {\overline{\boldsymbol{\beta}}}
{\longrightarrow}\Cl^{(t_{R})_{_{\!T}}}(T) \longrightarrow 0
$$
is split exact.
\end{coor}
\begin{proof} Recall that $(t_R)_\varphi=t_D$  \cite[Proposition 3.7]{fp}.  Then apply Theorem \ref{th:split}.
\end{proof}

    Note that,   when  we are dealing with a pullback
    diagram of type $(\square^{+})$,   $(t_{R})_{_{\!T}}\leq
    t_{T}$  (Lemma \ref{le:2.8} (2)) and so $\Cl^{(t_{R})_{_{\!T}}}(T)$ is a  subgroup of $\Cl^{t_{T}}(T)$.
In general,  it can happen that  $(t_{R})_{_{\!T}}\lneq t_{T}$ (for instance, when $M$ is not a $t_T$-ideal).
 We will show,   moreover, that $\Cl^{(t_{R})_{_{\!T}}}(T)$ can be a proper subgroup of $\Cl^{t_{T}}(T)$
 (Remark \ref{rk:2.20}).

\begin{coor} \label{cor:2.14} Under the same notation and hypotheses of
Corollary \ref{cor:2.13}, if we assume that $T$ is quasilocal, then
$\Cl^{(t_{R})_{_{\!T}}}(T)=0$. (In particular, we reobtain that
$\Cl^{t_{D}}(D) \cong \Cl^{t_{R}}(R)$,  see Corollary
\ref{cor:2.5}.)
\end{coor}
\begin{proof} Let $J$ be a ${(t_{R})_{_{\!T}}}$--invertible ${(t_{R})_{_{\!T}}}$--ideal
of $T$.  Then $J=(IT)^{{(t_{R})_{_{\!T}}}}$ for some nonzero finitely generated
fractional ideal $I$ of $R$ \cite[Proposition 2.6]{fpi}.  By the same
argument as in Claim 1 of the proof of Theorem \ref{thm1}, we have
$II^{-1}\not\subseteq M$.  Therefore $JJ^{-1}\supseteq
(IT)(IT)^{-1}=II^{-1}T=II^{-1}R_M=R_M=T$, and so $J$ is invertible in
$T$.  Since $T$ is quasilocal,  we conclude that $J$ is principal.
Therefore $\Cl^{{(t_{R})_{_{\!T}}}}(T)=0$.
\end{proof}

\begin{reem} \label{rk:2.20} \rm  
    %
 Note that for a pullback diagram of type $(\square^{+})$ with
$T$ quasilocal,  it is quite common that $\Cl^{t_T}(T)$ is nonzero,
but $\Ima(\overline{\boldsymbol{\beta}})=0$  (Corollaries
\ref{cor:2.13} and \ref{cor:2.14}).  An explicit example can be
obtained as follows.\   Let $T:= \mathbb{Q}[X^2, XY, Y^2]_{(X^2, XY,
Y^2)}$,  $M:= (X^2, XY, Y^2)T$, thus $T = \mathbb{Q}+M $, and  set
$R: = \mathbb{Z} +M$. Then, clearly $T =R_{M}$ and $M$ is a
$t_{R}$--prime of $R$. In this situation, the map
$\overline{\boldsymbol{\beta}}: \Cl^{t_{R}}(R) \rightarrow
\Cl^{t_{T}}(T) = \Cl^{t_{R_{M}}}(R_{M})$ is the zero map, while
$\Cl^{t_{T}}(T)$ is nonzero \cite[Proposition 2.3 and Example
3.4]{ar}.   Therefore in this case, by Corollary \ref{cor:2.14},
$\Cl^{(t_{R})_{_{\!T}}}(T)\neq \Cl^{t_{T}}(T)$.
\end{reem}

From Theorem \ref{th:split} applied to $\ast = d_{R}$, we reobtain
\cite[Theorem 2.5 (c)]{fg},  since
$(d_{R})_{\varphi} =d_{D}$ \cite[Proposition 3.3]{fp} and
$(d_{R})_{_{\!T}}=d_T$.   More precisely,

    \begin{coor}\label{pic}   Assume that we are dealing with a pullback
    diagram of type $(\square^{+})$.  Suppose that the map $\tilde{\varphi}:
    \mathcal U(T)\rightarrow k^\bullet/\mathcal U(D)$ is surjective.  Then
 $\Pic(R)\cong \Pic(D) \oplus \Pic(T)$.  \qed
\end{coor}

\begin{reem} \rm Note that, in \cite[Remark 2.7]{fg}, it was proved more
generally that:

\noindent \sl Assume that we are dealing with a pullback
    diagram of type $(\square)$.  The map $\tilde{\varphi}: \mathcal
    U(T)\rightarrow k^\bullet/\mathcal U(D)$ is surjective if and only if
    $\Pic(R)\cong \Pic(D) \oplus \Pic(T)$.  \rm

    \noindent  A similar result was reobtained in \cite[Theorem 3.9]{ac}.
   \end{reem}

   \smallskip

The next goal is to study  the behavior of  the property of being a
Pr\"{u}fer star multiplication domain in a pullback diagram of type
$(\square)$.  Recall that, given a star operation $\ast$ on an integral domain
$R$, we say that $R$ is a P$\ast$MD if for each nonzero
finitely generated  fractional ideal $I$ of $R$, $(II^{-1})^{\ast_{_{\!f}}} = R$ (cf.  for
instance \cite{fjs}, \cite{Gr}, \cite{k}, \cite{mz}, and \cite{hmm}).

\begin{thee} \label{th:2.22} Consider a pullback diagram of type $(\square)$ and let
$\ast$ be a star operation on $R$.  Then $R$ is a P$\ast$MD if and only
if $k$ is the quotient field of $D$, $D$ is a P${\ast}_{\varphi}$MD, $T$
is a P${({\ast})}_{_{\!T}}$MD, and $T_M$ is a valuation domain.
\end{thee}
\begin{proof} 
If $R$ is a P$\ast$MD, then $R$ is a
P$v$MD, and hence $k$ is the quotient field of $D$ and $T_M$ is a
valuation domain by \cite[Theorem 4.1]{fg}.  It is easy to see that if
$R$ is a P$\ast$MD, then $D$ is a P${\ast}_{\varphi}$MD and $T$ is a
P$({\ast})_{_{\!T}}$MD.
Actually, to prove that $T$ is a P$({\ast})_{_{\!T}}$MD,
let $J$ be a nonzero finitely generated ideal of $T$.
Since $T$ is $R$-flat, $J=IT$ for some finitely generated ideal $I$ of $R$.
Then by Lemma \ref{le:2.8} (4), $(JJ^{-1})^{((\ast)_{_{\!T}})_{_{\!f}}}
\supseteq (JJ^{-1})^{(\ast_{_{\!f}})_{_{\!T}}}
=(II^{-1}T)^{(\ast_{_{\!f}})_{_{\!T}}}=(II^{-1}T)^{\ast_{_{\!f}}}
=\left((II^{-1})^{\ast_{_{\!f}}}T\right)^{\ast_{_{\!f}}}=T^{\ast_{_{\!f}}}=T$.

 Conversely, assume that $k$ is the
quotient field of $D$, $D$ is a P${\ast}_{\varphi}$MD, $T$ is a
P${({\ast})}_{_{\!T}}$MD, and $T_M$ is a valuation domain. Since
$D$ and $T$ are P$v$MDs, $R$ is a P$v$MD by \cite[Theorem 4.1]{fg}.
Let $I$ be a nonzero finitely generated fractional ideal of $R$.
Then $(II^{-1})^{t_R}=R$, and hence $II^{-1}\not\subseteq M$. To
show that $I$ is  $\ast_{_{\!f}}$--invertible, we may assume that $I$ is a
nonzero finitely generated integral ideal of $R$ such that
$I\not\subseteq M$.  Since $D$ is a P$\ast_\varphi$MD,
$(\varphi(I)\varphi(I)^{-1})^{(\ast_\varphi)_{_{\!f}}}=D$.
Since $(\ast_\varphi)_{_{\!f}}=(\ast_{_{\!f}})_\varphi$ \cite[Proposition 3.6]{fp},
$(\varphi(I)\varphi(I)^{-1})^{(\ast_{_{\!f}})_\varphi}=D$, i.e.,
$\left((I+M)(I+M)^{-1}\right)^{\ast_{_{\!f}}}=R$, which implies that
$(II^{-1}+M)^{\ast_{_{\!f}}}=R$. Now suppose $II^{-1}\subseteq P$ for some
$P\in \Max^{\ast_{_{\!f}}}(R)$.  Then $M\not\subseteq P$, because otherwise
$R=(II^{-1}+M)^{\ast_{_{\!f}}}\subseteq P^{\ast_{_{\!f}}}=P$.
Note that $PT\in \Max^{((\ast)_{_{\!T}})_{_{\!f}}}(T)$ (by Claim 4 and 6
in the proof of Lemma \ref{le:2.8}).
But since $T$ is $R$-flat and $T$ is a P${({\ast})}_{_{\!T}}$MD,
$(IT(IT)^{-1})^{((\ast)_{_{\!T}})_{_{\!f}}}=(II^{-1}T)^{((\ast)_{_{\!T}})_{_{\!f}}}=T$, which
contradicts that $II^{-1}T\subseteq PT$.  Therefore
$II^{-1}\not\subseteq P$ for all $P\in \Max^{\ast_{_{\!f}}}(R)$, i.e.,
$(II^{-1})^{\ast_{_{\!f}}}=R$.  Thus $R$ is a P$\ast$MD.
\end{proof}

\begin{coor} \label{cor:2.24} Consider a pullback
diagram of type $(\square)$.  Then $R$ is a P$v_{R}$MD\ (= P$t_{R}$MD =
P$w_{R}$MD) if and only if $k$ is the quotient field of $D$, $D$ is a
P$v_{D}$MD \ (= P$t_{D}$MD =P$w_{D}$MD), $T$ is a P$(v_{R})_{_{\!T}}$MD \ (=
P$(t_{R})_{_{\!T}}$MD = P$(w_{R})_{_{\!T}}$MD), and $T_M$ is a valuation domain.
\end{coor}

\begin{proof} We can use Theorem \ref{th:2.22} and the following facts:
(1) for any star operation $\ast$ on an integral domain $A$,
$A$ is a P$\ast$MD
if and only if $A$ is a P$\tilde\ast$MD \cite[Theorem 3.1]{fjs};
(2) $(v_{R})_{\varphi}=v_{D}$ \cite[Corollary 2.13]{fp};
(3) when $k$ is the quotient field of $D$,
$\widetilde{( v_{R})_{_{\!T}}}= (\widetilde{\  v_{R}\ })_{_{\!T}} =  (w_{R})_{_{\!T}}
 \leq (t_{R})_{_{\!T}}= ((v_{R})_{_{\!f}})_{_{\!T}} \leq
 ((v_{R})_{_{\!T}})_{_{\!f}} \leq (v_{R})_{_{\!T}} $ (Lemma \ref{le:2.8}).
\end{proof}

\begin{reem} \label{rk:2.28} \rm Given a star operation $\ast$ on an integral domain
$R$, recall that  $R$  is a P$\ast$MD if and only if $R$ is a
P$v_{R}$MD and $\widetilde{\ast} =  t_{R}$ (or, equivalently,
${\ast}_{_{\!f}} =  t_{R}$) \cite[Proposition 3.4]{fjs}.  Therefore
(using Lemma \ref{le:2.8} (5) and \cite[Proposition 3.9]{fp})  the
previous theorem can be restated as follows: \sl  Consider a
pullback diagram of type $(\square)$ and let $\ast$ be a star
operation on $R$.  Then $ \widetilde{\ast}=  t_{R}$ and $R$ is a
P$v_{R}$MD if and only if $k$ is the quotient field of $D$,
${\widetilde{\ast}}_{\varphi}= t_{D}$,
${(\widetilde{\ast})}_{_{\!T}} = t_{T}$, $D$ is a P$v_{D}$MD, $T$ is
a P$v_T$MD, and $T_M$ is a valuation domain.
\end{reem}


\begin{leem}\label{le:2.29}
Let $R$ be a P$v_R$MD and let $T$ be a flat overring of $R$ such that $(R:T)\neq 0$.
Then $(w_{R})_{_{\!T}}= (t_{R})_{_{\!T}} = t_{T}=w_{T}$.
\end{leem}

 \begin{proof} Since $T$ is a flat overring of $R$, $T$ is a subintersection of $R$
and hence $T$ is a P$v_T$MD \cite[Theorem 3.11]{k}.
Recalling the fact that $w_A=t_A$ on a P$v_A$MD $A$ (\cite[Theorem 2.4]{p} or \cite[Proposition 3.4]{fjs}),
it suffices to show that $(t_{R})_{_{\!T}} = t_{T}$.

Note first that $T$ is a $w_R$-ideal of $R$ and hence a $t_R$-ideal of $R$.
Let $x\in T^{w_R}$. Then $xI \subseteq T$ for some finitely generated ideal
 $I$
of $R$ such that $I^{v_R}=R$  \cite[Remark 2.8]{FL2}.
By flatness,  $(IT)^{v_T}=(I^{v_R}T)^{v_T}=T$,   and thus $x\in T$.

Then $(t_{R})_{_{\!T}} \leq t_{T}$ and both are star operations on $T$ of finite type.
Let $J$ be a nonzero finitely generated integral ideal of $T$.
Then $J=IT$ for some finitely generated ideal $I$ of $R$.
By \cite[Proposition 2.17]{dhlz}, $I^{v_R}T$ is a $v_T$-ideal of $T$, and hence
$J^{t_T}=(I^{t_R}T)^{t_T}=(I^{v_R}T)^{t_T}=I^{v_R}T\subseteq (I^{v_R}T)^{t_R}=(I^{t_R}T)^{t_R}
=(IT)^{t_R}=J^{t_R}=J^{(t_R)_{_{\!T}}}$.
Thus we have $(t_{R})_{_{\!T}} = t_{T}$.
\end{proof}

\begin{coor}
Consider a pullback diagram of type $(\square)$.  Then $R$ is a
P$v_{R}$MD if and only if $k$ is the quotient field of $D$, $D$ is a
P$v_{D}$MD, $T$ is a P$v_{T}$MD, and $T_M$ is a valuation domain.
Moreover, in this situation, $(w_{R})_{_{\!T}}= (t_{R})_{_{\!T}} = t_{T}
=w_{T}$.
\end{coor}

 \begin{proof} The first statement is \cite[Theorem 4.1]{fg}
and the ``moreover'' statement follows from Lemma \ref{le:2.29}.
\end{proof}


\bigskip

\textsc{Dipartimento di Matematica,
 Universit\`a degli Studi ``Roma Tre'',
1, Largo San Leonardo Murialdo,
 00146 Roma, Italy}

 \it E-mail address: \rm \texttt{fontana@mat.uniroma3.it}

 \medskip

\textsc{Department of Mathematics, Chung-Ang University, Seoul
156-756, Korea}

\it  E-mail address: \rm  \texttt{mhpark@cau.ac.kr}

\end{document}